\documentclass[journal ]{new-aiaa}
\usepackage[utf8]{inputenc}
\usepackage{textcomp}

\usepackage{graphicx}
\usepackage{subcaption}
\usepackage{amsmath}
\usepackage{commath}
\usepackage{mathtools}
\usepackage[version=4]{mhchem}
\usepackage{siunitx}
\usepackage{longtable,tabularx}
\usepackage{cancel}
\usepackage{mathtools}

\usepackage{amsthm}
\newtheorem{remark}{Remark}
\newtheorem{assumption}{Assumption}
\newtheorem{problem}{Problem}
\newtheorem{lemma}{Lemma}

\newcommand{\comment}[1]{}

\title{On Covariance Steering with Nonlinear Measurements}

\author{Daniel C. Qi\footnote{Ph.D. Student, School of Aeronautics and Astronautics. (Corresponding Author: qi85@purdue.edu)} and Kenshiro Oguri\footnote{Assistant Professor, School of Aeronautics and Astronautics. Senior Member AIAA.}}
\affil{Purdue University, West Lafayette, Indiana, 47907}

\begin{document}

\maketitle

\begin{abstract}
This paper examines two ways nonlinear measurements complicate the statistical predictions during covariance steering. The covariance steering problem optimizes a control policy across a distribution of possible trajectories rather than a single specific trajectory. This approach has attracted growing interest as a way to account for uncertainty in aerospace applications such as trajectory optimization, stationkeeping, and landing. In many of these applications, covariance steering can account for measurement uncertainty through output feedback. Because standard formulations assume linear dynamics and measurements, applications with nonlinear measurements typically use an extended Kalman filter to estimate the state. This paper shows that techniques used to improve filter convergence, such as measurement underweighting, can reduce the accuracy of conventional covariance steering predictions by violating a key assumption. This is addressed with a new convex formulation that relaxes the assumption. This paper also demonstrates that adjusting the control policy during execution can improve agreement between predicted and actual state statistics. A numerical example of angles-only spacecraft rendezvous illustrates the benefits of the proposed method, and discussions are made regarding the implications of this generalized covariance steering formulation.
\end{abstract}

\section{Introduction}
This paper addresses discrete-time covariance steering with output-feedback under nonlinear measurements. Covariance steering (CS) is a stochastic optimal control problem that aims to control a distribution of trajectories, rather than a single trajectory. Although continuous-time formulations exist \cite{Hotz-CS, Chen-CS, Chen-CSOF}, finite-horizon discrete-time formulations are attractive to aerospace and astrodynamics applications that involve discrete maneuvers and measurement updates. The feasibility of this method has been demonstrated in problems such as trajectory optimization \cite{Ridderhof-SCP, Naoya-CS-JGCD, Jerry-CovSteering, Benedikter-CovSteering, Luca-CS, Gabriele-CS}, orbit stationkeeping \cite{Oguri-Safe-Autonomy}, proximity operations \cite{Fife-CS, Alkhenani-CS}, and landing \cite{Xiao-CS}. Covariance steering with output-feedback (CS-OF) \cite{Ridderhof-CSOF, Pilipovsky-CSOF} is a subset of the CS problem that considers control of a true distribution in the presence of measurement noise. This is relevant to many practical engineering applications in which the true state must be estimated from noisy measurements, a consideration included in many of the aforementioned examples \cite{Ridderhof-SCP, Naoya-CS-JGCD, Luca-CS, Gabriele-CS, Alkhenani-CS, Xiao-CS}.

Fundamentally, CS and CS-OF provide theoretical guarantees only for linear systems. In nonlinear systems, Gaussian and linearized approximations can produce reasonable results, but the approach can still break down in highly nonlinear systems if not treated properly \cite{Qi-min-NL}. Although methods tailored to nonlinear distribution steering and non-Gaussian control have been developed \cite{Qi-CUT-Steering, Qi-Risk-Measure, Burnett-NG-Control}, extending these formulations to output-feedback requires additional treatment of estimation uncertainty and its interaction with the control policy. Even with Gaussian measurement noise, the updated state distribution can be difficult to characterize analytically when the prior distribution is non-Gaussian. Thus, linearized CS-OF still remains attractive because of its computational tractability, motivating further research into mitigating nonlinear effects while preserving the computational advantages of this framework. One approach mitigates the effects of dynamical nonlinearities by minimizing the distribution's exposure to strongly nonlinear regions over the control horizon \cite{Qi-min-NL}. This helps maintain the accuracy of linearized uncertainty propagation. In CS-OF problems, nonlinearities can also arise in the measurement model, whose effects are to be studied in this paper.

Nonlinear measurements arise in many engineering applications. The extended Kalman filter (EKF) accommodates nonlinear measurement models through a first-order linearization of the system. This paper examines two aspects of the EKF in the context of CS-OF problems. Firstly, neglected higher-order terms can degrade estimation accuracy and, in severe cases, cause filter divergence. Modifications have therefore been developed to improve EKF robustness in nonlinear settings, namely measurement underweighting \cite{Zanetti-Underweighting}. Measurement underweighting deliberately modifies the optimal gain for the linearized model to reduce measurement influence and limit excessive contraction of the estimated error covariance in the presence of neglected measurement nonlinearities. This practice is extremely relevant, as its use can be seen on the Space Shuttle \cite{Shuttle-Underweight}, or with visual \cite{Park-Underweight} and angles-only \cite{Liao-Underweight} navigation. Second, the computed Kalman gain changes depending on the reference state about which the nonlinear measurement model is linearized. In an EKF, the measurement Jacobian is typically evaluated at a prior state estimate \cite{Stat-OD}. 

Although these considerations are well recognized in the field of navigation, the paper claims that effects arising from nonlinear measurements invalidate assumptions underlying existing CS-OF formulations, and therefore require explicit treatment within the covariance steering problem. The discrepancies associated with measurement underweighting and Kalman gain mismatch cannot be attributed solely to first-order approximation errors; addressing them requires modeling the state estimator and the feedback policy together in a single covariance steering formulation. This paper presents a method to solve this new CS-OF problem in a convex manner, and a modification to the control policy to account for discrepancies in Kalman gains. A spacecraft rendezvous example demonstrates how measurement underweighting can compromise covariance predictions in previous CS-OF formulations, and how the proposed framework with the correction term can resolve these errors. The paper concludes by examining the implications of the generalized CS-OF problem.

\section{Background}
The notation $I_{n}$ is an $n$-dimensional identity matrices and $0_{n\times m}$ is an $n\times m$ zero matrix. For brevity, the subscript in $0_{n\times m}$ is omitted whenever the dimensions of the surrounding matrices can be determined.

\subsection{Covariance Steering with Output-Feedback}\label{sec: CS-OF}
Let the time horizon be discretized into $N$ nodes, forming $N-1$ segments. Let $k = 0,1,\ldots, N-1$ denote the time index. The vectors $\boldsymbol{x}_k \in \mathbb{R}^{n_x}$, $\hat{\boldsymbol{x}}_k \in \mathbb{R}^{n_x}$, and $\hat{\boldsymbol{x}}^-_k \in \mathbb{R}^{n_x}$ denote the true state, filtered state estimate, and a priori state estimate at time step $t_k$, respectively. The problem of CS-OF treats each of these quantities as a random vector, and thus each has associated statistics. The mean of the true state vector at $k$ is denoted as $\boldsymbol{\mu}_k = \mathbb{E}[\boldsymbol{x}_k]$. The covariances of the true state and the filtered state at $k$ are then given by $P_k = \mathrm{Cov}(\boldsymbol{x}_k)$, $\hat{P}_k = \mathrm{Cov}(\hat{\boldsymbol{x}}_k)$, and $\hat{P}^-_k = \mathrm{Cov}(\hat{\boldsymbol{x}}^-_k)$. 

Consider a linear stochastic system with i.i.d. Gaussian white noise $\boldsymbol{w}_k \sim \mathcal{N}(0,I_{n_w})$ and control action $\boldsymbol{u}_k \in \mathbb{R}^{n_u}$. State dynamics, control mapping, and noise mapping matrices are given by $A_k\in \mathbb{R}^{n_x \times n_x}$, $B_k\in \mathbb{R}^{n_x \times n_u}$, and $G_k\in \mathbb{R}^{n_x \times n_w}$, respectively:
\begin{equation} \label{eq: true state dynamics}
\begin{aligned}
    \boldsymbol{x}_{k+1} &= A_k \boldsymbol{x}_k + B_k \boldsymbol{u}_k + G_k\boldsymbol{w}_k
\end{aligned}
\end{equation}
A linear measurement model for measurement $\boldsymbol{y}_k \in \mathbb{R}^{n_y}$ is given by a measurement mapping matrix $H_k \in \mathbb{R}^{n_y \times n_x}$ with i.i.d. Gaussian white noise $\boldsymbol{v}_k \sim \mathcal{N}(0, I_{n_v})$.
\begin{equation} \label{eq: linear msr model}
    \boldsymbol{y}_k = H_k \boldsymbol{x}_k + D_k \boldsymbol{v}_k 
\end{equation}

Let the estimation process be a linear update with a gain $L_k \in \mathbb{R}^{n_x \times n_y}$ for the a posteriori estimate in Eq.~\eqref{eq: filter update} and a time update for the a priori state estimate in Eq.~\eqref{eq: filter state dynamics}:
\begin{subequations}\label{eq: filter updates}
\begin{align}
\hat{\boldsymbol{x}}_k &= \hat{\boldsymbol{x}}^-_k + L_k(\boldsymbol{y}_k - H_k \hat{\boldsymbol{x}}^-_k) 
\label{eq: filter update} \\
\hat{\boldsymbol{x}}^-_{k+1} &= A_k\hat{\boldsymbol{x}}_k + B_k\boldsymbol{u}_k
\label{eq: filter state dynamics}
\end{align}
\end{subequations}

The boundary conditions of CS-OF are statistical in nature. The initial statistics of true and a priori estimate are normally distributed; both with mean $\boldsymbol{\mu}_0$, but separate covariances $P_{0}$ and $\hat{P}^-_{0}$ respectively. Terminal constraints are applied to constrain the final mean to be $\boldsymbol{\mu}_{f}$ and the final covariance to be less than $P_f$.
\begin{equation}\label{eq: inital and final conditions}
\begin{aligned}
\boldsymbol{x}_{0} \sim \mathcal{N}(\boldsymbol{\mu}_{0},P_{0})
,\qquad
\hat{\boldsymbol{x}}^-_{0} \sim \mathcal{N}(\boldsymbol{\mu}_{0},\hat{P}^-_{0}) 
,\qquad
\mathbb{E}[\boldsymbol{x}_{N-1}] = \boldsymbol{\mu}_{f}
,\qquad
\mathrm{Cov}(\boldsymbol{x}_{N-1}) \preceq P_f
\end{aligned} 
\end{equation}
For many applications, state-path and control magnitude constraints are needed. But since distributions such as Gaussian distributions can be unbounded, chance constraints are applied to ensure the desired constraint is satisfied with a specified high probability. Let $\mathcal{A}_\text{safe}$ be the safety set and $u_{\max}$ be the maximum allowed control action. The chance constraints are then
\begin{equation}\label{eq: chance constraints}
    \mathbb{P}\left(\boldsymbol{x}_k \in \mathcal{A}_\text{safe}\right) \geq 1 - \epsilon_x
    ,\qquad
    \mathbb{P}\left(\norm{\boldsymbol{u}_k}_2 \leq u_{\max} \right) \geq 1 - \epsilon_u
\end{equation}
where $\epsilon_x, \epsilon_u > 0$ are risk factors corresponding to state and control. Introduce an arbitrary objective function $J = \rho(\{\boldsymbol{x}_k,\boldsymbol{u}_k\})\in \mathbb{R}$ with $\rho$ mapping the the stochastic cost to $\mathbb{R}$. Problem~\ref{prob: stochastic CS-OF} states the general form of the CS-OF problem.
\begin{problem}[Covariance Steering with Output-Feedback]\label{prob: stochastic CS-OF}
Find $\{\boldsymbol{x}_k\}_{k=0,\ldots,N-1}$ and $\{\boldsymbol{u}_k\}_{k=0,\ldots,N-2}$ that minimizes $J$ subject to dynamics Eq.~\eqref{eq: true state dynamics}, measurement process Eq.~\eqref{eq: linear msr model}, filtering process Eq.~\eqref{eq: filter updates}, initial and terminal conditions Eq.~\eqref{eq: inital and final conditions}, and chance constraints Eq.~\eqref{eq: chance constraints}.
\end{problem}
Several assumptions are made to simplify the problem:
\begin{assumption}\label{assump: joint Gaussian}
$\boldsymbol{x}_0$ and $\hat{\boldsymbol{x}}^-_0$ are jointly Gaussian. 
\end{assumption}
This implies that the joint distribution exists. Also, it implies that their joint and marginal distributions are Gaussian and are uniquely determined by their mean and covariance \cite{Maybeck-StochasticTextbook}. In the context of CS-OF, controlling both the mean and covariance implies controlling the entire distribution.

\begin{assumption}\label{assump: filter first}
The estimator is designed prior to the control optimization problem. The estimation process ($L_k, H_k$) is not part of the control optimization process in Problem~\ref{prob: stochastic CS-OF}. 
\end{assumption}
A primary motivation for CS is the offline synthesis of feedback control policies \cite{Ridderhof-CSOF, Oguri-Safe-Autonomy, Naoya-CS-JGCD, Pilipovsky-CSOF}. Consistent with existing CS literature, Assumption~\ref{assump: filter first} treats the estimator structure and gains as prescribed inputs to the control optimization problem rather than as optimization variables. This is accomplished through a filtration process: given the initial filter error covariance $\widetilde{P}_0^-$ and any prescribed rule for computing $L_k$ from the prior covariance $\widetilde{P}_k^-$, the filter’s covariance prediction and update equations determine the gain histories used in Problem~\ref{prob: stochastic CS-OF}. The posterior covariance update is given by the Joseph's form: 
\begin{equation}\label{eq: Kalman Covariance Update}
 \widetilde{P}_k = (I_{n_x} - L_k H_k)\widetilde{P}^-_k (I_{n_x} - L_k H_k)^\top + L_k D_k D_k^\top L_k^\top
\end{equation}
where this covariance update holds for an arbitrary filter gain $L_k$. This is then followed by the time update:
\begin{equation}
    \widetilde{P}^-_{k+1} = A_k \widetilde{P}_k A_k^\top + G_k G_k^\top
\end{equation}
Unlike the aforementioned CS-OF studies, this assumption makes no assumptions about the optimality or functional form of $L_k$. This fact is examined more in depth in Section~\ref{sec: previous CSOF}.

\begin{assumption}\label{assump: unbiased estimate}
The estimator is unbiased. 
\end{assumption}
This assumption is reasonable, as most estimators are designed to be unbiased.

\subsubsection{Previous Works on CS-OF}\label{sec: previous CSOF}
There are many approaches to solve the CS and CS-OF problems. Methods differ in their covariance representation (full-covariance \cite{Pilipovsky-CSOF}, square-root covariance \cite{Naoya-Sqrt-CS}, or block square-root covariance \cite{Ridderhof-CSOF, Oguri-Safe-Autonomy}) and their choice of control policy. But in almost all existing CS-OF \cite{Chen-CSOF, Ridderhof-CSOF, Oguri-Safe-Autonomy, Naoya-CS-JGCD, Pilipovsky-CSOF, Xiao-CS}, an optimal Kalman gain is assumed:
\begin{equation}\label{eq: optimal kalman gain}
    L_{\text{kf},k} = \widetilde{P}^-_k H_k^\top(H_k\widetilde{P}^-_kH_k^\top + D_k D_k^\top)^{-1}
\end{equation}
The estimated state is then uncorrelated with the estimation error, which implies that the true state covariance is simply the sum of the estimated state covariance with Eq.~\eqref{eq: Kalman Covariance Update} \cite{Ridderhof-CSOF}:
\begin{equation}\label{eq: previous CSOF cov sum}
    P_k = \hat{P}_k + \widetilde{P}_k
\end{equation}
Under Assumption~\ref{assump: filter first}, $\widetilde{P}_k$ can be computed before the optimization problem. Then, previous works in CS-OF reduce the stochastic control problem to controlling $\hat{P}_k$, with an additive filter error process $\widetilde{P}_k$ when enforcing constraints on $P_k$. As mentioned earlier, there are different approaches to computing $\hat{P}_k$, but Eq.~\eqref{eq: previous CSOF cov sum} is the crucial commonality among existing CS-OF works. This assumption of uncorrelated estimates and estimation error is revisited in Section~\ref{sec: non-ortho CSOF}.

\subsection{Nonlinear Measurements}
Unlike Eq.~\eqref{eq: linear msr model}, many engineering examples feature nonlinear measurements. In these cases, consider the following nonlinear measurement model for a nonlinear measurement $\mathcal{Y}_k$:
\begin{equation}
    \mathcal{Y}_k = \mathcal{H}(\boldsymbol{x}_k) + D_k \boldsymbol{v}_k
\end{equation}
where $\boldsymbol{v}_k$ is still modeled as Gaussian additive noise, but involves a nonlinear measurement function $\mathcal{H}(\cdot)$. When applying the EKF, the measurement mapping matrix becomes the Jacobian of the nonlinear measurement function: 
\begin{equation}\label{eq: Hk linearization}
    H_k 
    =
    \frac{\partial}{\partial \boldsymbol{x}_k}
    \mathcal{H}(\boldsymbol{x}_k)
    \biggr\rvert_{\boldsymbol{x}_k = \boldsymbol{x}_k^*}
\end{equation}
where linearization is performed about a reference $\boldsymbol{x}_k^*$. EKFs can diverge in nonlinear systems due to this linearization inaccuracy; however, techniques can remedy certain diverging behaviors.

\subsubsection{Measurement Underweighting}\label{sec: msr uw}
One such technique is called measurement underweighting \cite{Zanetti-Underweighting}. A certain failure mode of EKFs is when a precise measurement is combined with large prior uncertainties. This combination can produce a large state correction based on a linearized measurement model. If the neglected higher-order terms are significant, the correction can be inaccurate and lead to filter divergence. Measurement underweighting mitigates this risk by reducing the influence of the measurement on the state estimate. While there are many ways to underweight a measurement \cite{Zanetti-Underweighting}, the underlying concept is to deliberately supply the filter estimate with a suboptimal Kalman gain. For this paper, consider the underweighted Kalman gain ${L}_{\mathrm{uw},k}$ in Eq.~\eqref{eq: uw Kalman gain}:
\begin{equation}\label{eq: uw Kalman gain}
    {L}_{\mathrm{uw},k} = \widetilde{P}^-_k H_k^\top
    \left( \frac{1}{p} H_k\widetilde{P}^-_kH_k^\top + D_k D_k^\top \right)^{-1}
\end{equation}
where $p \in (0,1]$ is the underweighting factor. As $p \rightarrow 1$, the underweight gain approaches the optimal Kalman gain, and as $p \rightarrow 0$, the measurement update is increasingly downweighted. With a suboptimal Kalman gain, only certain filter covariance updates, such as Eq.~\eqref{eq: Kalman Covariance Update}, can be used. 

\section{Key Consideration I: Non-Orthogonal CS-OF}\label{sec: non-ortho CSOF}
Recall Section~\ref{sec: previous CSOF} and specifically Eq.~\eqref{eq: previous CSOF cov sum}. Re-examining how the statistics of the true state are determined, the expression $\boldsymbol{x}_k =\hat{\boldsymbol{x}}_k+(\boldsymbol{x}_k - \hat{\boldsymbol{x}}_k)$ is true in any case. Taking the full expansion of the sum of two random variables yields:
\begin{equation}\label{eq: true state cov full sum}
    P_k = \hat{P}_k + \mathrm{Cov}(\boldsymbol{x}_k - \hat{\boldsymbol{x}}_k) 
    +    
    \mathrm{Cov}(\hat{\boldsymbol{x}}_k,\boldsymbol{x}_k - \hat{\boldsymbol{x}}_k) 
    +
    \mathrm{Cov}(\hat{\boldsymbol{x}}_k,\boldsymbol{x}_k - \hat{\boldsymbol{x}}_k)^\top
\end{equation}
where $\mathrm{Cov}(\hat{\boldsymbol{x}}_k,\boldsymbol{x}_k - \hat{\boldsymbol{x}}_k)  =  \mathbb{E}[\hat{\boldsymbol{x}}_k(\boldsymbol{x}_k - \hat{\boldsymbol{x}}_k)^\top] - \mathbb{E}[\hat{\boldsymbol{x}}_k] \mathbb{E}[(\boldsymbol{x}_k -\hat{\boldsymbol{x}}_k)^\top]$. From Assumption~\ref{assump: unbiased estimate}, $\mathbb{E}[(\boldsymbol{x}_k -\hat{\boldsymbol{x}}_k)^\top] = 0$, but the term $\mathbb{E}[\hat{\boldsymbol{x}}_k(\boldsymbol{x}_k - \hat{\boldsymbol{x}}_k)^\top]$ still remains. This term is only zero if the estimator satisfies the orthogonality principle, in which Eq.~\eqref{eq: true state cov full sum} reduces to Eq.~\eqref{eq: previous CSOF cov sum}. Since the traditional Kalman filter satisfies this principle, previous CS-OF need not consider these additional correlated effects. This paper refers to this case as \textit{orthogonal CS-OF}; that is, CS-OF methods that assume orthogonality between the estimate and estimate error. 

However, not all estimators satisfy the orthogonality principle. Section~\ref{sec: msr uw} presents a practical motivation for measurement underweighting, which employs a suboptimal Kalman gain and yields a nonzero cross-covariance between the state estimate and its estimation error. This motivates what this paper calls \textit{non-orthogonal CS-OF}; that is, CS-OF problems can have arbitrary statistical relationships between the estimate and the estimate error. In order to capture the arbitrary cross-correlation terms, an augmented state approach is considered.

\subsection{Augmented State Approach}
Define the augmented state vectors
\begin{equation}
    \mathbf{X}^-_k \triangleq 
    \begin{bmatrix}
        \boldsymbol{x}_k \\ \hat{\boldsymbol{x}}^-_k
    \end{bmatrix}
\qquad
    \mathbf{X}_k \triangleq 
    \begin{bmatrix}
        \boldsymbol{x}_k \\ \hat{\boldsymbol{x}}_k
    \end{bmatrix}
\end{equation}
For the covariance of the augmented state, define the following notation $\mathbf{P}^-_k = \mathrm{Cov}(\mathbf{X}^-_k)$ and $\mathbf{P}_k = \mathrm{Cov}(\mathbf{X}_k)$. In determining the evolution of the augmented system, the filter update from Eq.~\eqref{eq: filter update} can be written as
\begin{equation} \label{eq: augmented filter update}
\begin{aligned}
    \begin{bmatrix}
        \boldsymbol{x}_k \\ \hat{\boldsymbol{x}}_k
    \end{bmatrix}
    =
    \begin{bmatrix}
    I_{n_x} & 0 \\
    L_kH_k & I_{n_x} - L_k H_k
    \end{bmatrix}
    \begin{bmatrix}
        \boldsymbol{x}_k \\ \hat{\boldsymbol{x}}^-_k
    \end{bmatrix}
    +
    \begin{bmatrix}
        0 \\ L_k D_k
    \end{bmatrix}
    \boldsymbol{v}_k 
    \qquad \rightarrow \qquad
    \mathbf{X}_k = \Lambda_k \mathbf{X}^-_k 
    +
    \mathbf{L}_k
    \boldsymbol{v}_k 
\end{aligned}
\end{equation}
The block matrices $\Lambda_k$ and $\mathbf{L}_k$ are introduced to simplify notation. Likewise, the dynamics and control update from Eqs.~\eqref{eq: true state dynamics} and \eqref{eq: filter state dynamics} is written as
\begin{equation} \label{eq: augmented control update}
\begin{aligned}
\begin{bmatrix}
    \boldsymbol{x}_{k+1} \\ \hat{\boldsymbol{x}}^-_{k+1}
\end{bmatrix}
=
    \begin{bmatrix}
    A_k & 0 \\
    0 & A_k 
    \end{bmatrix}
    \begin{bmatrix}
        \boldsymbol{x}_k \\ \hat{\boldsymbol{x}}_k
    \end{bmatrix}
    +
    \begin{bmatrix}
        B_k \\ B_k
    \end{bmatrix}
    \boldsymbol{u}_k
        +    
    \begin{bmatrix}
        G_k \\ 0
    \end{bmatrix}\boldsymbol{w}_k 
\qquad \rightarrow \qquad
    \mathbf{X}^-_{k+1} = 
    \widetilde{A}_k \mathbf{X}_k + \widetilde{B}_k\boldsymbol{u}_k 
    +    
    \widetilde{G}_k
    \boldsymbol{w}_k 
\end{aligned}
\end{equation}
with similar block matrices $\widetilde{A}_k$, $\widetilde{B}_k$, and $\widetilde{G}_k$. Eqs.~\eqref{eq: augmented filter update} and \eqref{eq: augmented control update} can then be combined into one equation to determine the evolution of $\mathbf{X}_{k}$ to $\mathbf{X}_{k+1}$:
\begin{equation}\label{eq: augmented state evolution combined}
\begin{aligned}
    \mathbf{X}_{k+1} &= 
    \mathbf{A}_k \mathbf{X}_{k}
    +
    \mathbf{B}_k\boldsymbol{u}_k
    +    
    \mathbf{G}_k
    \boldsymbol{w}_k 
    +
    \mathbf{L}_{k+1}
    \boldsymbol{v}_{k+1}
\end{aligned}
\end{equation}
where $\mathbf{A}_k = \Lambda_{k+1}\widetilde{A}_k$, $\mathbf{B}_k = \Lambda_{k+1}\widetilde{B}_k$, and $\mathbf{G}_k = \Lambda_{k+1}\widetilde{G}_k$. These block matrices can be viewed as the matrices that capture the aggregate effects of dynamics and controls followed by the filter update. The next step is to determine the evolution of the mean and covariance for this augmented system.

\subsection{Mean Evolution}
If an unbiased estimator is used (Assumption~\ref{assump: unbiased estimate}), then $\mathbb{E}[\boldsymbol{x}_k] = \mathbb{E}[\hat{\boldsymbol{x}}^-_k] = \mathbb{E}[\hat{\boldsymbol{x}}_k]$ for all $k$. Then the mean evolution for the augmented updates can be easily computed through linearity of expectation, and in fact Eqs.~\eqref{eq: augmented filter update} and \eqref{eq: augmented control update} can be decoupled to just track the mean of the true state:
\begin{equation}\label{eq: mean evolution}
    \boldsymbol{\mu}_{k+1} = A_k \boldsymbol{\mu}_{k} + B_k \bar{\boldsymbol{u}}_k
\end{equation}
where $\mathbb{E}[\boldsymbol{u}_k] = \bar{\boldsymbol{u}}_k$. Using this, define an aggregated mean vectors $\bar{\mathbf{X}}^-_k(=\mathbb{E}[\mathbf{X}^-_k])=\bar{\mathbf{X}}_k(=\mathbb{E}[\mathbf{X}_k])= [\boldsymbol{\mu}_k^\top,\boldsymbol{\mu}_k^\top]^\top$. While computing the mean evolution is straightforward, the crux of most CS algorithms is how the control policy influences the covariance evolution. The following section presents a method where the evolution of the covariance is convex with respect to the control gains with an appropriately chosen control policy.

\subsection{Convex Solution to Covariance Evolution: Block Matrix with State-History Feedback}
A block matrix approach has been applied to many previous CS problems \cite{Ridderhof-CSOF, Oguri-Safe-Autonomy, Okamoto-CS, Okamoto-CS-Vehicle}. This paper applies a similar approach as \cite{Ridderhof-CSOF} in its treatment of the block matrices to retrieve a convex expression with respect to a control policy. Firstly, Eq.~\eqref{eq: augmented state evolution combined} can be expressed in the block matrix form:
\begin{equation}
\begin{aligned}
    \begin{bmatrix}
        \mathbf{X}_0 \\
        \mathbf{X}_1 \\
        \mathbf{X}_2 \\
        \vdots
    \end{bmatrix}
    =&
    \begin{bmatrix}
        \Lambda_0 \\
        \mathbf{A}_0\Lambda_0 \\
        \mathbf{A}_1\mathbf{A}_0\Lambda_0 \\
        \vdots
    \end{bmatrix}
    \mathbf{X}_0^-
    +
    \begin{bmatrix}
        0 & 0 & 0 & \cdots \\
        \mathbf{B}_0 & 0 & 0 & \cdots \\
        \mathbf{A}_1\mathbf{B}_0 & \mathbf{B}_1 & 0 & \cdots \\
        \vdots & \vdots & \vdots & \ddots
    \end{bmatrix}
    \begin{bmatrix}
        \boldsymbol{u}_0 \\
        \boldsymbol{u}_1 \\
        \boldsymbol{u}_2 \\
        \vdots
    \end{bmatrix}
    +
    \\
    &
    \begin{bmatrix}
        0 & 0 & 0 & \cdots \\
        \mathbf{G}_0 & 0 & 0 & \cdots \\
        \mathbf{A}_1\mathbf{G}_0 & \mathbf{G}_1 & 0 & \cdots \\
        \vdots & \vdots & \vdots & \ddots
    \end{bmatrix}
    \begin{bmatrix}
        \boldsymbol{w}_0 \\
        \boldsymbol{w}_1 \\
        \boldsymbol{w}_2 \\
        \vdots
    \end{bmatrix}
    +
    \begin{bmatrix}
        \mathbf{L}_0
            & 0 & 0 & 0 & \cdots \\
            \mathbf{A}_0\mathbf{L}_0 & \mathbf{L}_1 & 0 & 0 & \cdots \\
           \mathbf{A}_1\mathbf{A}_0\mathbf{L}_0
            & \mathbf{A}_1\mathbf{L}_1
            & \mathbf{L}_2 & 0 & \cdots \\
        \vdots & \vdots & \vdots & \vdots & \ddots
    \end{bmatrix}
    \begin{bmatrix}
        \boldsymbol{v}_0 \\
        \boldsymbol{v}_1 \\
        \boldsymbol{v}_2 \\
        \vdots
    \end{bmatrix}
\end{aligned}
\end{equation}
Written in compact form:
\begin{equation}\label{eq: block matrix original}
    \mathbf{X} = \mathbf{A} \mathbf{X}_0^- + \mathbf{B}\mathbf{u} + \mathbf{G} \mathbf{w}+ \mathbf{L} \mathbf{v}
\end{equation}
with stacked vectors $\mathbf{X} = [\mathbf{X}_0^\top, \ldots, \mathbf{X}_{N-1}^\top]^\top$, $\mathbf{u} = [\boldsymbol{u}_0^\top, \ldots, \boldsymbol{u}_{N-2}^\top]^\top$, and similarly for $\mathbf{w}$ and $\mathbf{v}$. Define the control policy to be a linear state-history feedback (SHF) policy, and express it in terms of the augmented state:
\begin{equation}\label{eq: control policy}
    \boldsymbol{u}_k = \bar{\boldsymbol{u}}_k + \sum^k_{i = 0} K_{k,i}(\hat{\boldsymbol{x}}_i - \boldsymbol{\mu}_i)
    =
    \bar{\boldsymbol{u}}_k + \sum^k_{i = 0} 
    \begin{bmatrix}
    0_{n_u \times n_x} &  K_{k,i}
    \end{bmatrix}
    \left(
    \mathbf{X}_i -
    \bar{\mathbf{X}}_i
    \right)
\end{equation}
Then, the corresponding block matrix form for control can be given by
\begin{equation}\label{eq: block matrix control}
    \mathbf{u} = \bar{\mathbf{u}} + \mathbf{K}(\mathbf{X} - \bar{\mathbf{X}})
\end{equation}
where $\bar{\mathbf{u}} = [\bar{\boldsymbol{u}}_0^\top, \ldots, \bar{\boldsymbol{u}}_{N-2}^\top]^\top$, $\bar{\mathbf{X}} = [\bar{\mathbf{X}}_0^\top,\ldots, \bar{\mathbf{X}}_{N-1}^\top]^\top$, and
\begin{equation}\label{eq: sparse K}
    \mathbf{K} =
    \begin{bmatrix}
        0 & K_{0,0} & 0 & 0
        & \cdots & 0 & 0 & 0 & 0 \\
        0 & K_{1,0} & 0 & K_{1,1}
        & \cdots & 0 & 0 & 0 & 0 \\
        \vdots & \vdots & \vdots & \vdots
        & \ddots & \vdots & \vdots & \vdots & \vdots \\
        0 & K_{N-2,0} & 0 & K_{N-2,1}
        & \cdots & 0 & K_{N-2,N-2} & 0 & 0
    \end{bmatrix}
\end{equation}
Combining Eq.~\eqref{eq: block matrix original} with Eq.~\eqref{eq: block matrix control} yields
\begin{equation}
    \mathbf{X} = \mathbf{A} \mathbf{X}_0^- + \mathbf{B}\bar{\mathbf{u}} + \mathbf{B}\mathbf{K}(\mathbf{X} - \bar{\mathbf{X}})+ \mathbf{G} \mathbf{w}+ \mathbf{L} \mathbf{v}
\end{equation}
By linearity of expectation, the mean of the stacked vector is $\bar{\mathbf{X}} = \mathbf{A} \bar{\mathbf{X}}_0^- + \mathbf{B}\bar{\mathbf{u}}$. Standardizing gives 
\begin{equation}
    \mathbf{X} - \bar{\mathbf{X}}
    =
    \mathbf{A} (\mathbf{X}_0^- - \bar{\mathbf{X}}_0^- ) + \mathbf{B}\mathbf{K}(\mathbf{X} - \bar{\mathbf{X}}) + \mathbf{G} \mathbf{w}+ \mathbf{L} \mathbf{v}
\end{equation}
Since $\mathbf{B}$ is strictly lower triangular and $\mathbf{K}$ is block lower triangular, $(I - \mathbf{B}\mathbf{K})$ is invertible. Then, term $\mathbf{X} - \bar{\mathbf{X}}$ can be isolated to the right-hand side:
\begin{equation}\label{eq: X minus Xbar, nonconvex}
    \mathbf{X} - \bar{\mathbf{X}}
    =
    (I - \mathbf{B}\mathbf{K})^{-1}
    [\mathbf{A} (\mathbf{X}_0^- - \bar{\mathbf{X}}_0^- ) + \mathbf{G} \mathbf{w}+ \mathbf{L} \mathbf{v}]
\end{equation}
Define a new decision variable $\mathbf{F} = \mathbf{K}(I-\mathbf{B}\mathbf{K})^{-1}$ \cite{Skaf-affine-controller, Ridderhof-CSOF} that satisfies 
\begin{equation}\label{eq: K to F sub}
    I+\mathbf{B}\mathbf{F} = (I-\mathbf{B}\mathbf{K})^{-1}
\end{equation}
The gain can be uniquely recovered with $\mathbf{K} = \mathbf{F}(I+\mathbf{B}\mathbf{F})^{-1}$. Lemma~\ref{lemma: block F} in Appendix~\ref{appendix: sparse F} shows that $\mathbf{F}$ has the same sparse structure as $\mathbf{K}$ in Eq.~\eqref{eq: sparse K}, so the approach is then to optimize over $\mathbf{F}$ to remove the nonconvexity from the inverse. Substituting Eq.~\eqref{eq: K to F sub} into Eq.~\eqref{eq: X minus Xbar, nonconvex} yields
\begin{equation}\label{eq: X minus Xbar, convex}
    \mathbf{X} - \bar{\mathbf{X}}
    =
    (I + \mathbf{B}\mathbf{F})
    [\mathbf{A} (\mathbf{X}_0^- - \bar{\mathbf{X}}_0^- ) + \mathbf{G} \mathbf{w}+ \mathbf{L} \mathbf{v}]
\end{equation}
Since $\mathbf{w}$ and $\mathbf{v}$ are both stacked vectors of i.i.d. Gaussian white noise, the stacked noise vectors themselves are also independent with identity covariance. Define $\mathbf{S}$ to be the covariance of $\mathbf{A} (\mathbf{X}_0^- - \bar{\mathbf{X}}_0^- ) + \mathbf{G} \mathbf{w}+ \mathbf{L} \mathbf{v}$ from Eq.~\eqref{eq: X minus Xbar, convex}:
\begin{equation}
    \mathbf{S} = \mathrm{Cov}( \mathbf{A} (\mathbf{X}_0^- - \bar{\mathbf{X}}_0^- )  + \mathbf{G} \mathbf{w}+ \mathbf{L} \mathbf{v}) = \mathbf{A} \mathbf{P}_0^{-}\mathbf{A}^\top + \mathbf{G}\mathbf{G}^\top + \mathbf{L}\mathbf{L}^\top
\end{equation}
It follows that the covariance of the stacked augmented vector $\mathbf{P} = \mathrm{Cov}(\mathbf{X})$ can be expressed as a square-root covariance, and likewise for the stacked control covariance $\mathbf{P}_u = \mathrm{Cov}(\mathbf{u})$
\begin{equation}\label{eq: cov evolution}
    \mathbf{P}^{1/2} = (I + \mathbf{B}\mathbf{F})\mathbf{S}^{1/2}
    ,\qquad
    \mathbf{P}_u^{1/2} = \mathbf{F}\mathbf{S}^{1/2}
\end{equation}
where $\mathbf{P} = \mathbf{P}^{1/2}(\mathbf{P}^{1/2})^\top$, $\mathbf{P}_u = \mathbf{P}_u^{1/2}(\mathbf{P}_u^{1/2})^\top$, and $\mathbf{S} = \mathbf{S}^{1/2}(\mathbf{S}^{1/2})^\top$.
Thus, it can be seen that both the history of the true state and control covariances, captured by $\mathbf{P}^{1/2}$ and $\mathbf{P}_u^{1/2}$, are affine in $\mathbf{F}$. To retrieve the covariances at a specific $t_k$, define the selection matrix $E_x(k)$ that selects the $k$-th index of the true state from the stack augmented state ($\boldsymbol{x}_k = E_x(k)\mathbf{X}$). Likewise, define the selection matrix for control ($\boldsymbol{u}_k = E_u(k)\mathbf{u}$). With this, the covariance for the true state $P_k$ and the control $P_{u_k} = \mathrm{Cov}(\boldsymbol{u}_k)$ at $t_k$ can be determined in square-root form $P^{1/2}_k = E_x(k) \mathbf{P}^{1/2}$ and $P^{1/2}_{u_k} = E_u(k) \mathbf{P}_u^{1/2}$.

\subsection{Chance Constraints}\label{sec: chance constraints}
Chance constraints are generally nonconvex and are often replaced by conservative deterministic approximations. These approximations may remain nonconvex, depending on the CS formulation \cite{Pilipovsky-CSOF, Naoya-CS-JGCD, Rapakoulias-CS-CC}. However, the block square-root covariance framework allows certain chance constraints to be expressed in convex form \cite{Oguri-Safe-Autonomy, Ridderhof-CSOF}. This section illustrates their incorporation into the proposed formulation by focusing on those used in the numerical example in Section~\ref{sec: numerical example}. Define the function $m_{\chi^2}(\epsilon,n)$ as the quantile function of a chi-square distribution with $n$ degrees of freedom evaluated at the probability value of $1-\epsilon$. Similarly, define $m_{\mathcal{N}}(\epsilon)$ as the quantile function of the standard normal distribution at the probability value of $1-\epsilon$. The matrix $H_r$ selects the 3-dimensional position vector from the state. 

\subsubsection{Control Magnitude}
The control magnitude constraint is a 2-norm constraint on the stochastic variable: $\mathbb{P}\left(\norm{\boldsymbol{u}_k}_2 \leq u_{\max} \right) \geq 1 - \epsilon_u$. Expressed using the mean and covariance \cite{Oguri-Safe-Autonomy}:
\begin{equation}\label{eq: cc, control mag}
   \norm{\bar{\boldsymbol{u}}_k}_2 + \sqrt{m_{\chi^2}(\epsilon_u,n_u)}\norm{P_{u_k}^{1/2}}_2 \leq u_{\max}
\end{equation}

\subsubsection{Approach Angle Constraint}
An approach angle constraint can be considered as a cone path constraint on position \cite{Oguri-Safe-Autonomy}:
\begin{equation}
    \mathbb{P}\left(
    \norm{A_{\text{cone}}H_r\boldsymbol{x}_k}_2
    -
    \boldsymbol{b}_{\text{cone}}^\top H_r\boldsymbol{x}_k
    \leq
    0
    \right) \geq 1 - \epsilon_x
\end{equation}
where $A_{\text{cone}}$ and $\boldsymbol{b}_{\text{cone}}$ defines the geometry of the cone. The numerical example in Section~\ref{sec: numerical example} considers a conical approach region with its apex at
the origin and its symmetry axis aligned with the $-y$-axis.
Given these specification, the cone geometry is defined as $A_{\mathrm{cone}} =
\begin{bsmallmatrix}
1 & 0 & 0 \\
0 & 0 & 1
\end{bsmallmatrix}$ and
$b_{\mathrm{cone}} = [0,-\tan\theta_{\mathrm{cone}},0]^\top$ for a cone half-angle of $\theta_{\text{cone}}$. Expressing the approach angle constraint using the mean and covariance \cite{Oguri-Safe-Autonomy}:
\begin{equation}\label{eq: cc, cone}
   \norm{A_{\text{cone}}H_r\boldsymbol{\mu}_k}_2 
   - \boldsymbol{b}^\top_{\text{cone}}H_r\boldsymbol{\mu}_k
  + \sqrt{m_{\chi^2}\left(\frac{\epsilon_x}{2},2\right)}\norm{A_{\text{cone}}H_rP^{1/2}_k}_2
  + m_{\mathcal{N}}\left(\frac{\epsilon_x}{2}\right)\norm{\boldsymbol{b}^\top_{\text{cone}}H_r P^{1/2}_k}_2 \leq 0
\end{equation}

\subsubsection{Covariance Constraint}
Although not a chance constraint, a constraint on the final covariance is generally desirable for CS-OF problems. In a square-root form this is given by \cite{Oguri-Safe-Autonomy, Okamoto-CS-Vehicle, Ridderhof-CSOF} 
\begin{equation}\label{eq: final cov constraint}
    P_{N-1} \preceq P_f
    \qquad\rightarrow \qquad
    \norm{P_f^{-1/2} P^{1/2}_{N-1}}_2 \leq 1
\end{equation}

\subsubsection{Delta-V 99}
A common objective problem within astrodynamics is called $\Delta V_{99}$, which is the 99th percentile of $\Delta V$ costs. Computing the quantile is nonconvex, so an upper-bounding $\Delta V_{99,\mathrm{ub}}$ is used in place \cite{Oguri-Safe-Autonomy}:
\begin{equation}\label{eq: deltaV 99}
\Delta V_{99,\mathrm{ub}}=\sum_{k=0}^{N-2}
   \norm{\bar{\boldsymbol{u}}_k}_2 + \sqrt{m_{\chi^2}(0.01,n_u)}\norm{P_{u_k}^{1/2}}_2
\end{equation}

\section{Key Consideration II: Reference Kalman Gain v.s. Realized Kalman Gain}
CS-OF yields two distinctions: a \emph{linear} solution and a \emph{linearized} solution. A linear solution holds exactly for linear stochastic systems with linear measurements under the assumptions in Section~\ref{sec: CS-OF}. For nonlinear systems, the dynamics and/or the measurement functions are linearized to obtain a linearized CS-OF solution. For example, the measurement Jacobian from Eq.~\eqref{eq: Hk linearization} is evaluated about the reference mean during CS-OF optimization \cite{Naoya-CS-JGCD, Oguri-Safe-Autonomy}. This approximates the evolution of the state statistics using Jacobians evaluated along a reference trajectory. For brevity, this paper refers to these Jacobians as the reference matrices.

For linear dynamics and measurements, the Kalman gain is independent of the a priori state estimate when the a priori covariance and measurement noise covariance are fixed. In an extended Kalman filter, however, the measurement Jacobian is evaluated at $\hat{\boldsymbol{x}}_k^-$, so different realizations of the a priori state estimates can yield different realized Kalman gains, even for the same a priori covariance. Because conventional linearized CS-OF computes control actions using the reference matrices, this mismatch can lead to discrepancies between the predicted and realized statistics. This paper addresses these discrepancies by introducing a gain correction into the control policy.

\subsection{Gain Correction Term}
Let $L_k^*$ and $H_k^*$ be the reference Kalman gain and measurement Jacobian used in the linearized CS-OF optimization, and $\boldsymbol{\hat{x}}^*_k$ and ${\boldsymbol{\hat{x}}_k^-}^*$ be the estimates produced with those reference matrices. Specifically, $\boldsymbol{\hat{x}}^*_k = {\boldsymbol{\hat{x}}_k^-}^* + L_k^*[\boldsymbol{y}_k - \mathcal{H}({\boldsymbol{\hat{x}}_k^-}^*)]$. Similarly, let $L_k$ and $H_k$ be the realized Kalman gain and measurement Jacobian obtained by the nonlinear estimator, and $\boldsymbol{\hat{x}}_k$ and ${\boldsymbol{\hat{x}}_k^-}$ be the estimates produced with the realized matrices. Let $\Tilde{\boldsymbol{y}}_k$ and $\Tilde{\boldsymbol{y}}^*_k$ be the realized and reference innovations respectively. Introduce an auxiliary variable $\boldsymbol{z}_k$ that accounts for this mismatch such that $\boldsymbol{z}_k =\boldsymbol{\hat{x}}_k-\boldsymbol{\hat{x}}^*_k$ and $\boldsymbol{z}^-_k =\boldsymbol{\hat{x}}^-_k-{\boldsymbol{\hat{x}}_k^-}^*$. Then the evolution of this auxiliary variable is
\begin{subequations}
\begin{align}
    \boldsymbol{z}_k &=
    \boldsymbol{z}^-_{k} + L_k \Tilde{\boldsymbol{y}}_k - L_k^* \Tilde{\boldsymbol{y}}_k^* 
    \label{eq: zk update, nonlinear H}
    \\
    \boldsymbol{z}^-_{k+1} &= A_k \boldsymbol{z}_k
\end{align}
\end{subequations}
where $\boldsymbol{z}^-_{0} =0$. The innovation are computed as $\Tilde{\boldsymbol{y}}_k^* = \boldsymbol{y}_k - \mathcal{H}({\boldsymbol{\hat{x}}_k^-}^*)$ and $\Tilde{\boldsymbol{y}}_k = \boldsymbol{y}_k - \mathcal{H}({\boldsymbol{\hat{x}}_k^-})$, or $\Tilde{\boldsymbol{y}}_k^* - \Tilde{\boldsymbol{y}}_k =\mathcal{H}({\boldsymbol{\hat{x}}_k^-}) - \mathcal{H}({\boldsymbol{\hat{x}}_k^-}^*)$. Under first-order reference linearization, $\Tilde{\boldsymbol{y}}_k^* - \Tilde{\boldsymbol{y}}_k \approx H_k^* \boldsymbol{z}^-_k$. Substituting into Eq.~\eqref{eq: zk update, nonlinear H} yields
\begin{equation}\label{eq: zk update, linear H}
    \boldsymbol{z}_k \approx
    (I-L_k^*H_k^*)\boldsymbol{z}^-_{k} + (L_k - L_k^*) \Tilde{\boldsymbol{y}}_k
\end{equation}
The choice of using Eq.~\eqref{eq: zk update, nonlinear H} or \eqref{eq: zk update, linear H} is up to the user, and both account for the difference between the realized Kalman gain $L_k$ and the reference Kalman gain $L_k^*$. In Eq.~\eqref{eq: zk update, nonlinear H}, the correction to $\boldsymbol{z}_k$ uses measurement predictions computed with the nonlinear function $\mathcal{H}(\cdot)$. For highly nonlinear measurements, higher-order terms in these predictions can amplify deviations in $\boldsymbol{z}_k$ and, through the feedback policy, in the control input. Eq.~\eqref{eq: zk update, linear H} therefore replaces the nonlinear measurement term in the $\boldsymbol{z}_k$ update with its first-order approximation to reduce the gain correction's sensitivity to large deviations. This paper uses Eq.~\eqref{eq: zk update, linear H} as the $\boldsymbol{z}_k$ update.

Nevertheless, the control law optimized with linearized CS-OF is based on $\boldsymbol{\hat{x}}^*_k$ (such as the SHF policy in Eq.~\eqref{eq: control policy} when applied to nonlinear measurements, for example). Thus, this correction term is applied to the SHF policy:
\begin{equation}
    \boldsymbol{u}_k = \bar{\boldsymbol{u}}_k + \sum^k_{i = 0} K_{k,i}(\hat{\boldsymbol{x}}^*_i- \boldsymbol{\mu}_i)=\bar{\boldsymbol{u}}_k + \sum^k_{i = 0} K_{k,i}(\hat{\boldsymbol{x}}_i - \boldsymbol{z}_i - \boldsymbol{\mu}_i)
\end{equation}
The application to the SHF policy is only used as an example, but it can be seen that the evolution of $\boldsymbol{z}_k$ is independent of the control policy. As a result, substitution of $\boldsymbol{\hat{x}}^*_k = \hat{\boldsymbol{x}}_k - \boldsymbol{z}_k$ can be performed for any control policy that uses $\boldsymbol{\hat{x}}^*_k$, for example the Markovian feedback policy $\boldsymbol{u}_k = \bar{\boldsymbol{u}}_k + K_k (\hat{\boldsymbol{x}}^*_k- \boldsymbol{\mu}_k)$ from other CS-OF methods \cite{Naoya-CS-JGCD, Pilipovsky-CSOF}.

This correction is applied during closed-loop operations to mitigate the effects of gain mismatch arising from nonlinear measurements without altering either the orthogonal or non-orthogonal CS-OF optimization. However, because the EKF ultimately relies on first-order linearizations, exact agreement between the predicted and realized statistics is not guaranteed.

\section{Numerical Example: Angles-Only Spacecraft Rendezvous}\label{sec: numerical example}
The example uses the 6-dimensional linear Clohessy-Wiltshire (CWH) equation as the dynamics. The state is $\boldsymbol{x}_k = [x_k,y_k,z_k,\dot{x}_k,\dot{y}_k,\dot{z}_k]^\top$ with the origin at the chief spacecraft. Impulsive control is used. The linear dynamics and control mapping are 
\begin{equation}
A_k = e^{A \Delta t_k}
,\qquad
B_k = A_k B
\end{equation}
where $\Delta t_k = t_{k+1}-t_k$. The matrix $A$ is the continuous-time equations of motion for linear CWH, and $B$ is the impulsive control mapping matrix. These matrices are 
\begin{equation}
A =
\begin{bmatrix}
0_{3\times3}& I_3\\
A_{21} & A_{22}
\end{bmatrix}
,\qquad 
A_{21}  =
\begin{bmatrix}
3\tilde{n}^2 & 0 & 0 \\
0    & 0 & 0 \\
0    & 0 & -\tilde{n}^2
\end{bmatrix}
,\qquad 
A_{22}  =
\begin{bmatrix}
0   & 2\tilde{n} & 0 \\
-2\tilde{n} & 0  & 0 \\
0    & 0 & 0
\end{bmatrix}
,\qquad 
B
=
\begin{bmatrix}
    0_{3\times 3}\\
    I_3
\end{bmatrix}
\end{equation}
where $\tilde{n}$ is the mean motion of the circular chief orbit. Let the nonlinear measurement be the azimuth $\theta \in [-\pi,\pi)$ and elevation $\phi\in [-\pi/2,\pi/2]$ angles of the \emph{chief observed by the deputy}. The measurement model is
\begin{equation}
\begin{aligned}
    \begin{bmatrix}
        \theta \\ \phi
    \end{bmatrix} = \mathcal{H}_{(\theta,\phi)}(\boldsymbol{x}_k) + D_k \boldsymbol{v}_k,
    \qquad
   \mathcal{H}_{(\theta,\phi)}(\boldsymbol{x}_k) = 
    \begin{bmatrix}
        \mathrm{atan2}(-y,-x)\\
        \mathrm{atan2}\left(-z,\sqrt{x^2+y^2}\right)
    \end{bmatrix} 
\end{aligned}
\end{equation}
where $\mathrm{atan2}(\cdot,\cdot)$ is the four-quadrant inverse tangent function. For CS-OF optimization, this measurement Jacobian is evaluated about the reference mean $H_k^* = \frac{\partial}{\partial \boldsymbol{x}_k}\mathcal{H}_{(\theta,\phi)}(\boldsymbol{x}_k) \biggr\rvert_{\boldsymbol{x}_k = \boldsymbol{\mu}_k^*}$ like prior works \cite{Oguri-Safe-Autonomy, Naoya-CS-JGCD}, but evaluated about the realized a priori estimate $H_k = \frac{\partial}{\partial \boldsymbol{x}_k}\mathcal{H}_{(\theta,\phi)}(\boldsymbol{x}_k) \biggr\rvert_{\boldsymbol{x}_k = \boldsymbol{x}^-_k}$ during each realized trial.

\begin{remark}
The numerical example uses linear dynamics to isolate the effect of measurement nonlinearity on the accuracy of CS-OF predictions.
\end{remark}

The numerical example assumes that \textit{initially} the statistics of the estimate and estimation error are orthogonal. The formulation also allows non-orthogonal initial statistics, as Section~\ref{sec: discussion} explains. For the orthogonal case, Appendix~\ref{appendix: initial cov} derives the following initial augmented covariance:
\begin{equation}\label{eq: numerical example, initial cov}
    \mathbf{P}^-_0 = 
\begin{bmatrix}
    \hat{P}^-_0 + \widetilde{P}^-_0 &\hat{P}^-_0\\
    \hat{P}^-_0 & \hat{P}^-_0
\end{bmatrix}
\end{equation}
The chance constraints in Section~\ref{sec: chance constraints} are imposed, and the $\Delta V_{99,\mathrm{ub}}$ is used as the objective function. 

The rendezvous example considers a V-bar approach to the chief, with prescribed waypoints along the V-bar. The deputy's mean position must pass through these waypoints, and control actions are permitted only at the waypoint nodes and the final control node $N-2$. These requirements yield the following constraints:
\begin{equation}\label{eq: contraints vbar}
\begin{aligned}
    \bar{\boldsymbol{u}}_k &= 0,
    \qquad E_u(k)\mathbf{F} = 0,
    \qquad
    && \forall k \in
    \{0,\ldots,N-2\}
    \setminus
    \left(\mathbb{Z}_{\text{V-bar}}\cup\{N-2\}\right)
    \\
    H_r\boldsymbol{\mu}_{k_i}
    &= [0, \bar{y}_i,0]^{\top},
    && i=1,\ldots,M
\end{aligned}
\end{equation}
where $\mathbb{Z}_{\text{V-bar}}=\{k_1,\ldots,k_M\}$ is the set of prescribed waypoint node indices, $\bar{y}_i$ is the along-track coordinate of the $i$-th waypoint. Imposing $\bar{\boldsymbol{u}}_k=0$ eliminates the nominal control, while $E_u(k)\mathbf{F}=0$ implies $E_u(k)\mathbf{K}=0$ and therefore eliminates the feedback contribution. Together, these constraints ensure that $\boldsymbol{u}_k=0$ at every node outside $\mathbb{Z}_{\text{V-bar}}\cup\{N-2\}$. Measurements are still acquired at each optimization node.

\begin{problem}[Convex Non-Orthogonal CS-OF for V-bar Rendezvous Example]
Find $\{\boldsymbol{\mu}_k\}_{k=0,\ldots,N-1}$, $\{\bar{\boldsymbol{u}}_k\}_{k=0,\ldots,N-2}$, $\mathbf{P}^{1/2},\mathbf{P}_u^{1/2}, \mathbf{F}$ that minimizes Eq.~\eqref{eq: deltaV 99} subject to mean dynamics Eq.~\eqref{eq: mean evolution}, covariance dynamics Eq.~\eqref{eq: cov evolution}, control magnitude change constraint Eq.~\eqref{eq: cc, control mag}, approach angle chance constraint Eq.~\eqref{eq: cc, cone}, V-bar assignments Eq.~\eqref{eq: contraints vbar}, initial and terminal mean constraints $\boldsymbol{\mu}_{0} = \bar{\boldsymbol{\mu}}_{0}$, $\boldsymbol{\mu}_{N-1} = \boldsymbol{\mu}_{f}$, and initial and terminal covariance conditions Eqs.~\eqref{eq: numerical example, initial cov}, \eqref{eq: final cov constraint}
\end{problem}

The rendezvous is conducted around Earth with a gravitational parameter of $\mu_{\text{Earth}} \approx 398600$ km$^3$/s$^2$. The uncertainty in the angle measurement is meant to be a realistic $\sigma_{\angle} = 20$ arcsec \cite{Kruger-angles}, leading to $D_k = \text{diag}([\sigma_{\angle},\sigma_{\angle}])$. The number of optimization nodes is $N=14$ evenly spaced in time with $\Delta t = 60$ sec. The rest of the system parameters are given in
Table~\ref{tab: system Parameters} and constraint values in Table~\ref{tab: Constraint Parameters}. All length units are nondimensionalized by $0.1$ km and time units by $100$ sec for better numerical conditioning in the optimization.
\begin{table}[htbp]
	\fontsize{10}{10}\selectfont
    \caption{System Parameters}
   \label{tab: system Parameters}
        \centering 
   \begin{tabular}{l  c  c} 
      \hline 
      Parameter & Value & Units\\
      \hline 
      Semi-major Axis of Chief Orbit & 7228 & km \\
      Initial Mean ($\bar{\boldsymbol{\mu}}_{0}$) & $[0,-1.5,0,0,0,0]^\top$ & [km, km/s] \\  
      Initial Estimate Covariance ($\hat{P}^-_0$) & $\text{diag}([30,30,30,10,10,10]^2)$ & [m, cm/s]$^2$ \\
      Initial Filter Error Covariance ($\widetilde{P}^-_0$) & $\text{diag}([50,30,50,1,1,1]^2)$ & [m, cm/s]$^2$ \\
      Dynamical Stochastic Process ($G_k$) & $\text{diag}([1,1,1,1,1,1])$ & [cm, mm/s] \\
    \hline
   \end{tabular}
\end{table}

\begin{table}[htbp]
	\fontsize{10}{10}\selectfont
    \caption{Constraint Parameters}
   \label{tab: Constraint Parameters}
        \centering 
   \begin{tabular}{l  c  c} 
      \hline 
      Parameter & Value & Units\\
      \hline 
      Target Mean ($\boldsymbol{\mu}_{f}$) & $[0,-0.20,0,0,0,0]^\top$ & [km, km/s] \\
      Terminal Covariance ($P_f$) & $\text{diag}([15,15,15,10,10,10]^2)$ & [m, cm/s]$^2$ \\
      Chance Constraint Risk ($\epsilon_x,\epsilon_u$) & $10^{-3}$ & - \\
      Max Magnitude Control ($u_{\max}$) & 10 & m/s \\
      Cone Half-Angle ($\theta_{\text{cone}}$) & 10 & degrees \\
      V-bar Nodes ($\mathbb{Z}_{\text{V-bar}}$) & $\{0,4,9\}$& -\\
      V-bar Targets ($\bar{y}_i$) & $\{-1.5,-1,-0.5\}$& km \\
    \hline
   \end{tabular}
\end{table}

\begin{figure}[htb]
	\centering\includegraphics[width=0.7\textwidth]{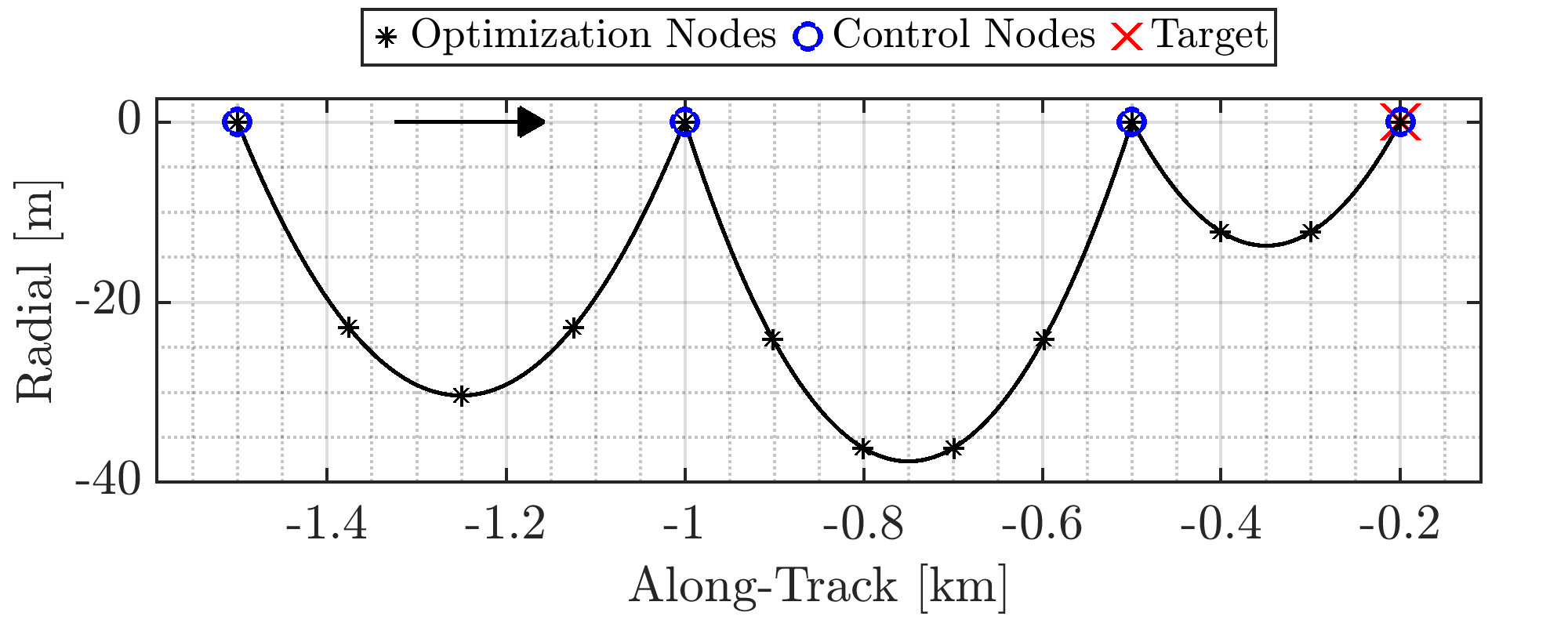}
    \caption{Reference trajectory of V-bar approach \emph{(axes different scales)}.}
	\label{fig: traj}
\end{figure}

The linearized CS-OF formulation requires an initial reference trajectory to evaluate the measurement Jacobians. The deterministic trajectory computed using the mean values, shown in Figure~\ref{fig: traj}, serves as this reference. Because the CS-OF solution can shift the mean trajectory when accounting for the covariance, the reference and associated Jacobians may require updates within a sequential convex framework. However, for the system and constraint parameters considered here, a single convex solve sufficed for each CS-OF case because the optimized mean trajectory matched the initial reference.

\subsection{Methodology}
To motivate the importance of this paper's considerations from nonlinear measurements, the traditional orthogonal frameworks are applied for comparison. This analysis applies the control policies found from (1) \emph{orthogonal} CS-OF \emph{without} underweighting and (2) \emph{orthogonal} CS-OF \emph{with} underweighting. Orthogonal CS-OF with measurement underweighting uses the conventional relation $P_k = \hat{P}_k + \widetilde{P}_k$, hence naively neglecting the cross-covariance terms. The estimation error covariance $\widetilde{P}_k$ is computed using Eq.~\eqref{eq: Kalman Covariance Update} with the underweighted Kalman gain ${L}_{\mathrm{uw},k}$. The orthogonal CS-OF algorithm chosen for comparison is from Ridderhof et al. \cite{Ridderhof-CSOF} as it also uses the same state-history feedback policy in Eq.~\eqref{eq: control policy}. Both the orthogonal and non-orthogonal problems possess identical objective functions, initial/final conditions, and constraints. To add, conventional orthogonal CS-OF formulations do not consider gain correction from nonlinear measurement. Accordingly, this correction is omitted from the orthogonal CS-OF approach evaluated here unless explicitly stated otherwise.

\subsection{Need for Underweighting}\label{sec: results, ortho CSOF}
\begin{figure}[!htb]
  \centering 
    \subcaptionbox
    {Positional Distribution at Terminal Time.\label{fig: ortho_CSOF, dist}}
    {\includegraphics[width=0.48\textwidth]{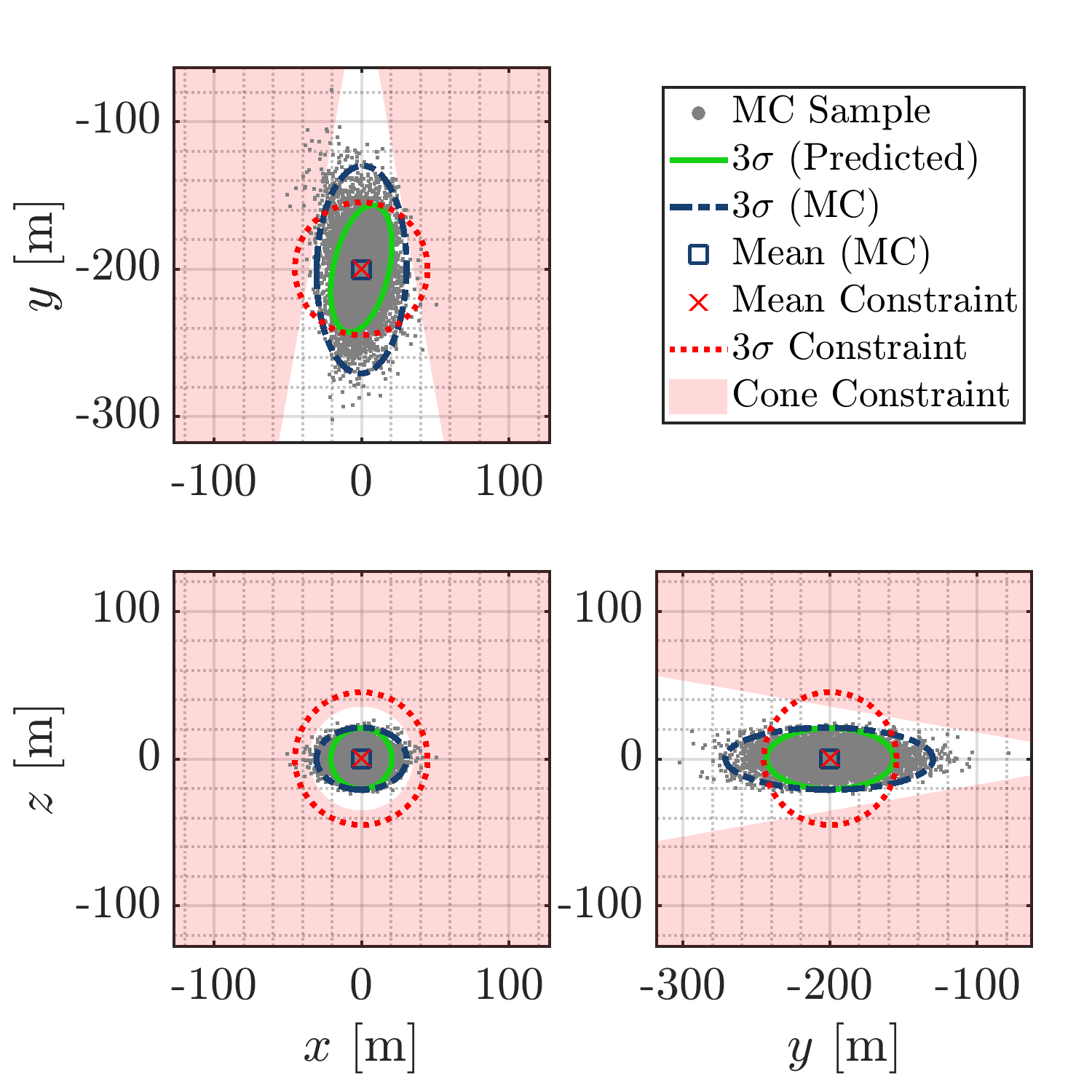}}
    \hskip 0.1truein
    \subcaptionbox
    {Time History of Filter Statistics.\label{fig: ortho_CSOF, filter}}
    {\includegraphics[width=0.48\textwidth]{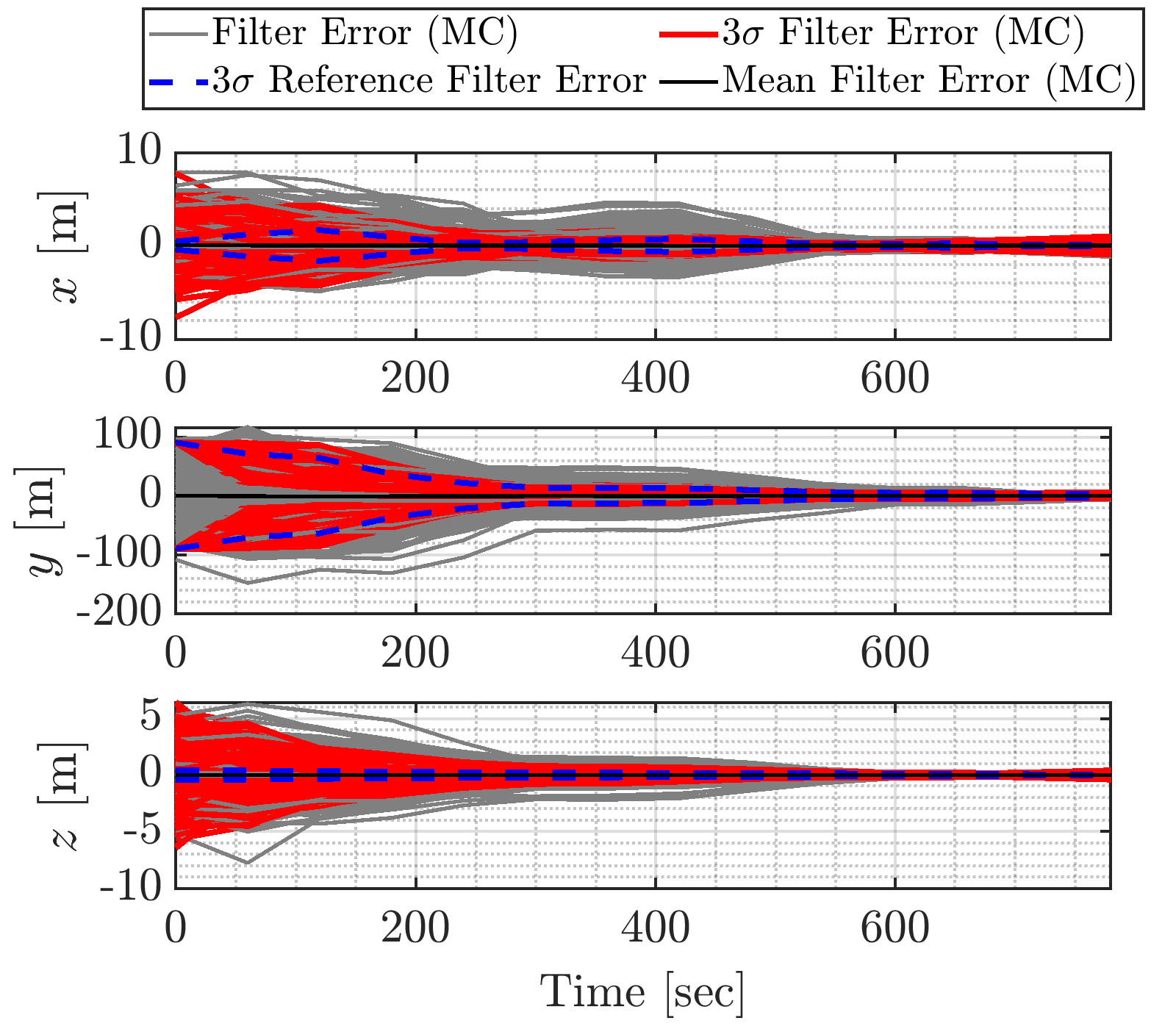}}
    \caption{Monte Carlo ($n_{\text{samples}}=10,000$): Nonlinear simulation results using orthogonal CS-OF without measurement underweighting.}
    \label{fig: ortho_CSOF}
\end{figure}

Figure~\ref{fig: ortho_CSOF} \footnote{For visual clarity, only the first 1,000 samples are shown in the filter statistics plots. All samples are included in the position distribution plots.} shows the resulting distribution using the control gains computed from orthogonal CS-OF with no underweighting. It is clear from Figure~\ref{fig: ortho_CSOF, dist} that the predicted covariances do not align with the actual covariance from the Monte Carlo. The terminal covariance constraint is violated, with a significant number of samples also violating the cone constraint. The reasoning behind this inaccuracy can be inferred from Figure~\ref{fig: ortho_CSOF, filter}: while the estimate itself appears to be unbiased, the reference estimation error covariance used to compute the covariance of the true state differs from many of the realized filter error covariances. Furthermore, many of the samples exceed the filter's predicted $3\sigma$ error interval at the beginning, leading to further statistical inconsistencies in the estimation error. While the errors ultimately converge in the steady state, the initial divergence is operationally undesirable. From Figure~\ref{fig: ortho_CSOF, filter}, it can be concluded that measurement underweighting is needed.   

\begin{figure}[!htb]
  \centering 
    \subcaptionbox
    {Positional Distribution at Terminal Time.\label{fig: ortho_CSOF_uw, dist}}
    {\includegraphics[width=0.48\textwidth]{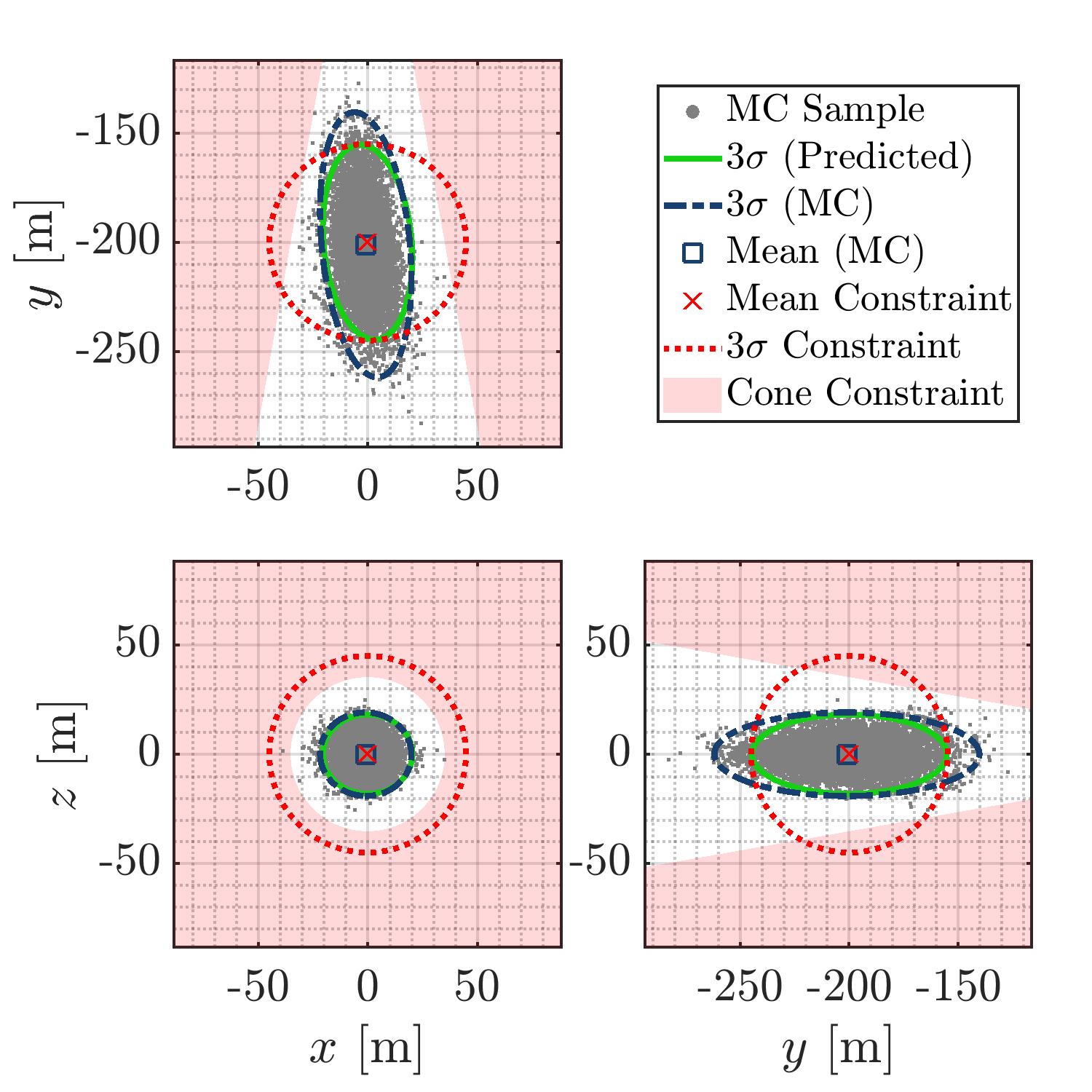}}
    \hskip 0.1truein
    \subcaptionbox
    {Time History of Filter Statistics.\label{fig: ortho_CSOF_uw, filter}}
    {\includegraphics[width=0.48\textwidth]{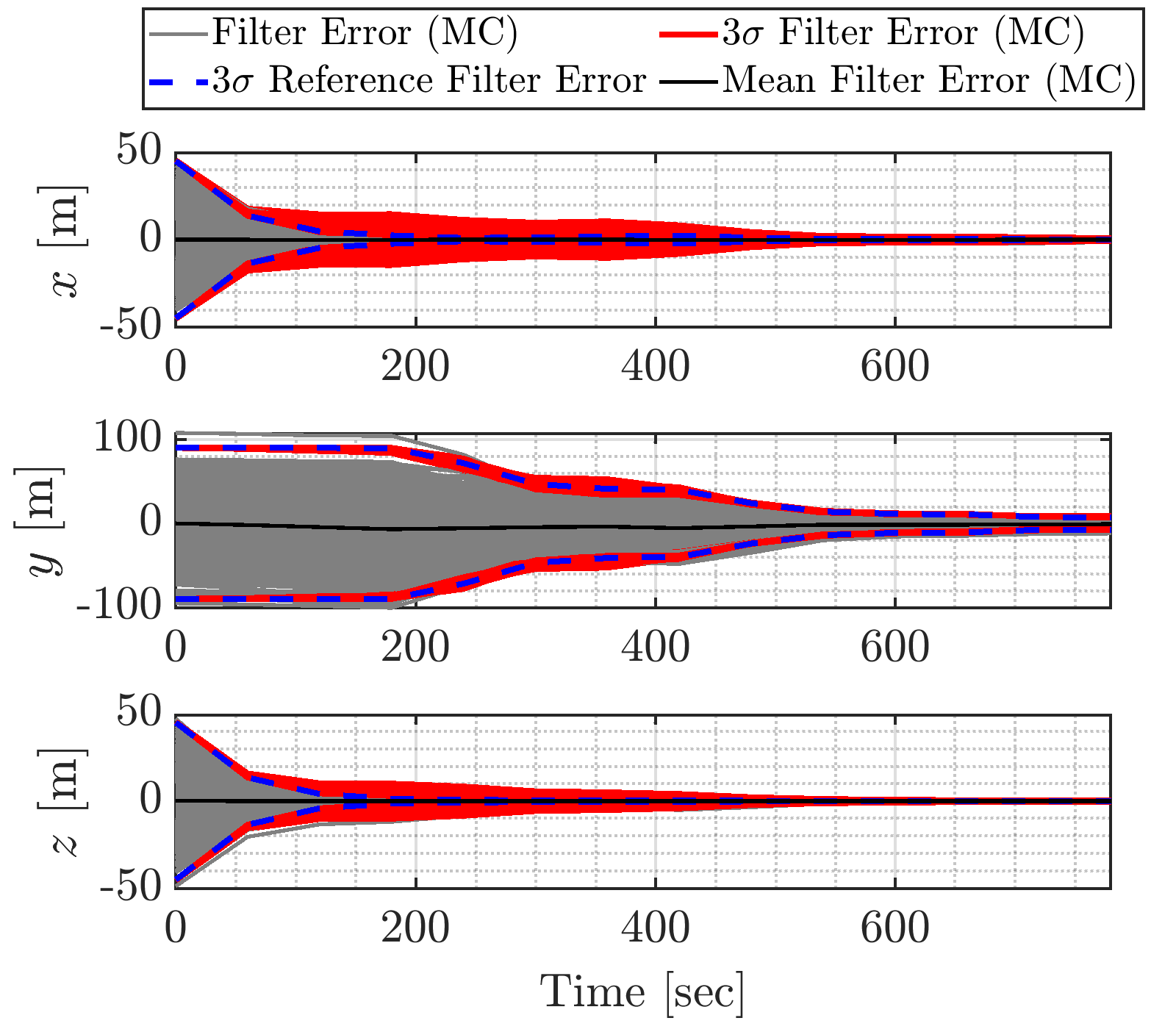}}
    \caption{Monte Carlo ($n_{\text{samples}}=10,000$): Nonlinear simulation results using orthogonal CS-OF with measurement underweighting.}
    \label{fig: ortho_CSOF_uw}
\end{figure}

Figure~\ref{fig: ortho_CSOF_uw} shows the resulting distribution using the control gains computed from orthogonal CS-OF with an underweighting factor of $p=0.7$, which is in the neighborhood of operational values \cite{Shuttle-Underweight, Zanetti-Underweighting}. The introduction of measurement underweighting substantially improved the performance of the filter, as shown in Figure~\ref{fig: ortho_CSOF_uw, filter}. Compared with the previous Figure~\ref{fig: ortho_CSOF, filter}, the reference filter error covariance is more consistent with the realized filter error covariances, and the $3\sigma$ intervals of the unweighted case effectively capture the realized errors. However, there is still a discrepancy found in Figure~\ref{fig: ortho_CSOF_uw, dist}: the actual covariance does not satisfy the terminal covariance constraint, as seen in the $x$-$y$ and $y$-$z$ joint marginal distributions. It can be concluded that while the filter in the underweighted case is more statistically consistent, the orthogonal CS-OF solution does not capture the additional cross-correlation produced by the suboptimality.

\subsection{Need for Non-Orthogonal CS-OF}
\begin{figure}[!htb]
  \centering 
    \subcaptionbox
    {Positional Distribution at Terminal Time.\label{fig: nonortho_CSOF, dist}}
    {\includegraphics[width=0.48\textwidth]{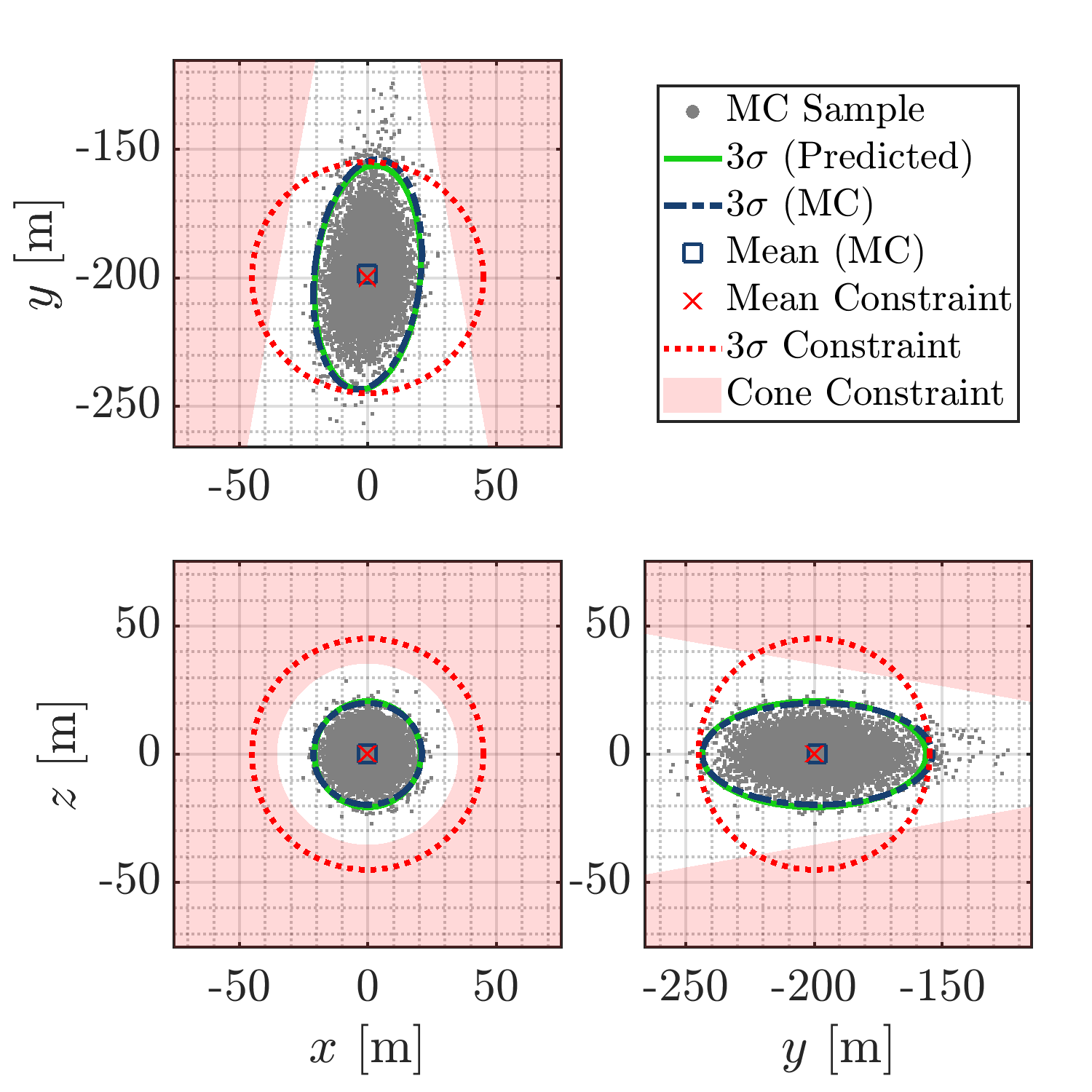}}
    \hskip 0.1truein
    \subcaptionbox
    {Time History of Filter Statistics.\label{fig: nonortho_CSOF, filter}}
    {\includegraphics[width=0.48\textwidth]{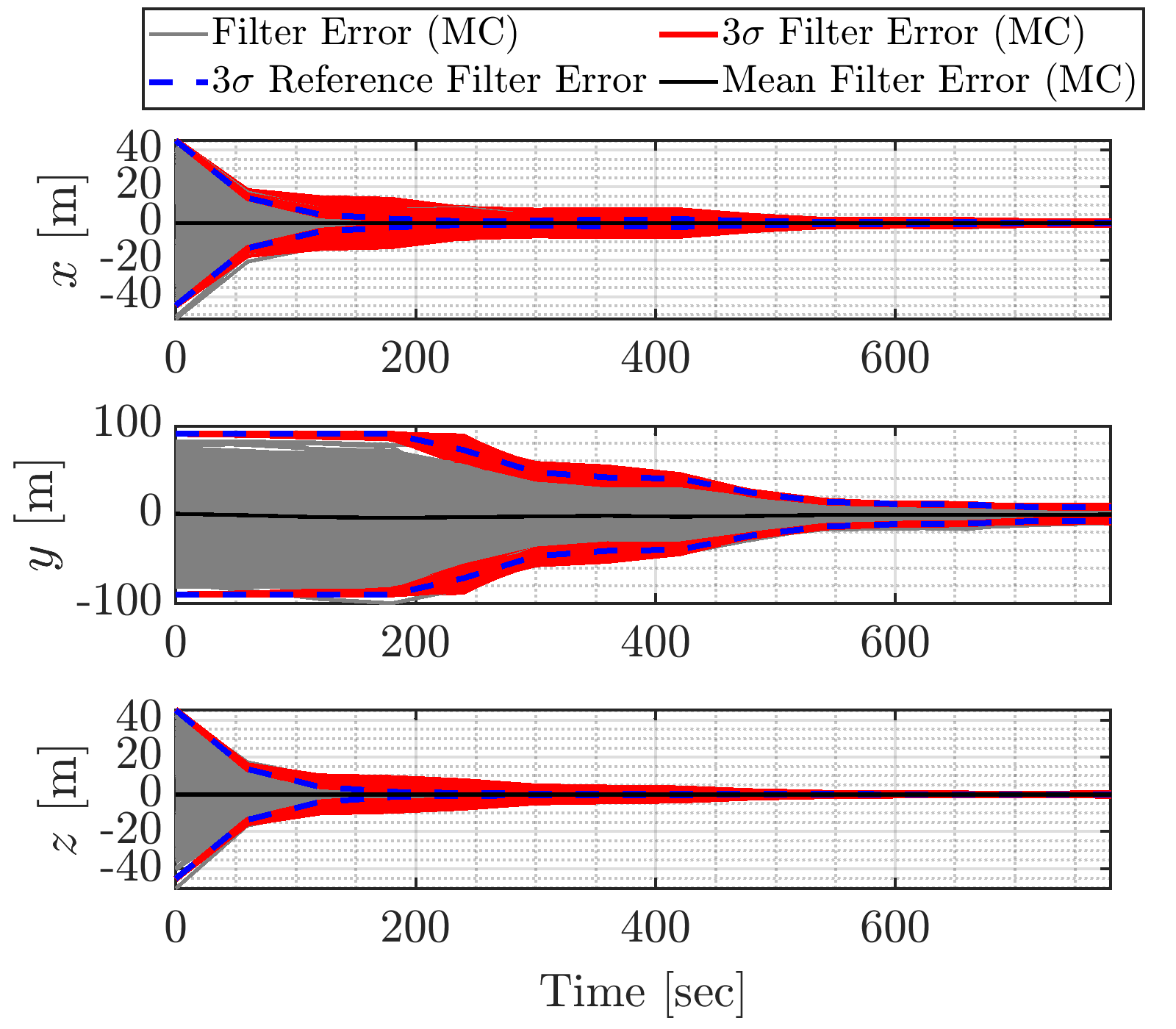}}
    \caption{Monte Carlo ($n_{\text{samples}}=10,000$): Nonlinear simulation results using non-orthogonal CS-OF with measurement underweighting and gain correction.}
    \label{fig: nonortho_CSOF}
\end{figure}

Figure~\ref{fig: nonortho_CSOF} shows the resulting distribution using the control gains computed from non-orthogonal CS-OF along with measurement underweighting and the gain correction within the control policy. The underweighting results in a similarly converged filter as the underweighted orthogonal CS-OF case, but the main difference is in the statistical consistency in Figure~\ref{fig: nonortho_CSOF, dist}. Recall from the underweighted orthogonal CS-OF case, Figure~\ref{fig: ortho_CSOF_uw, dist} shows that the actual covariance violated the terminal covariance constraint. But after applying the non-orthogonal CS-OF formulation, the actual covariance is consistent with the predicted covariance, and the violation in the terminal covariance constraint is resolved. This highlights the stronger coupling between covariance control and filtering in nonlinear systems. Techniques that improve filter consistency can also affect covariance control, so these interactions must be accounted for to obtain a statistically consistent solution.

\begin{figure}[!htb]
  \centering 
    \subcaptionbox
    {Boxchart of $\Delta V$ History.\label{fig: control history}}
    {\includegraphics[width=0.48\textwidth]{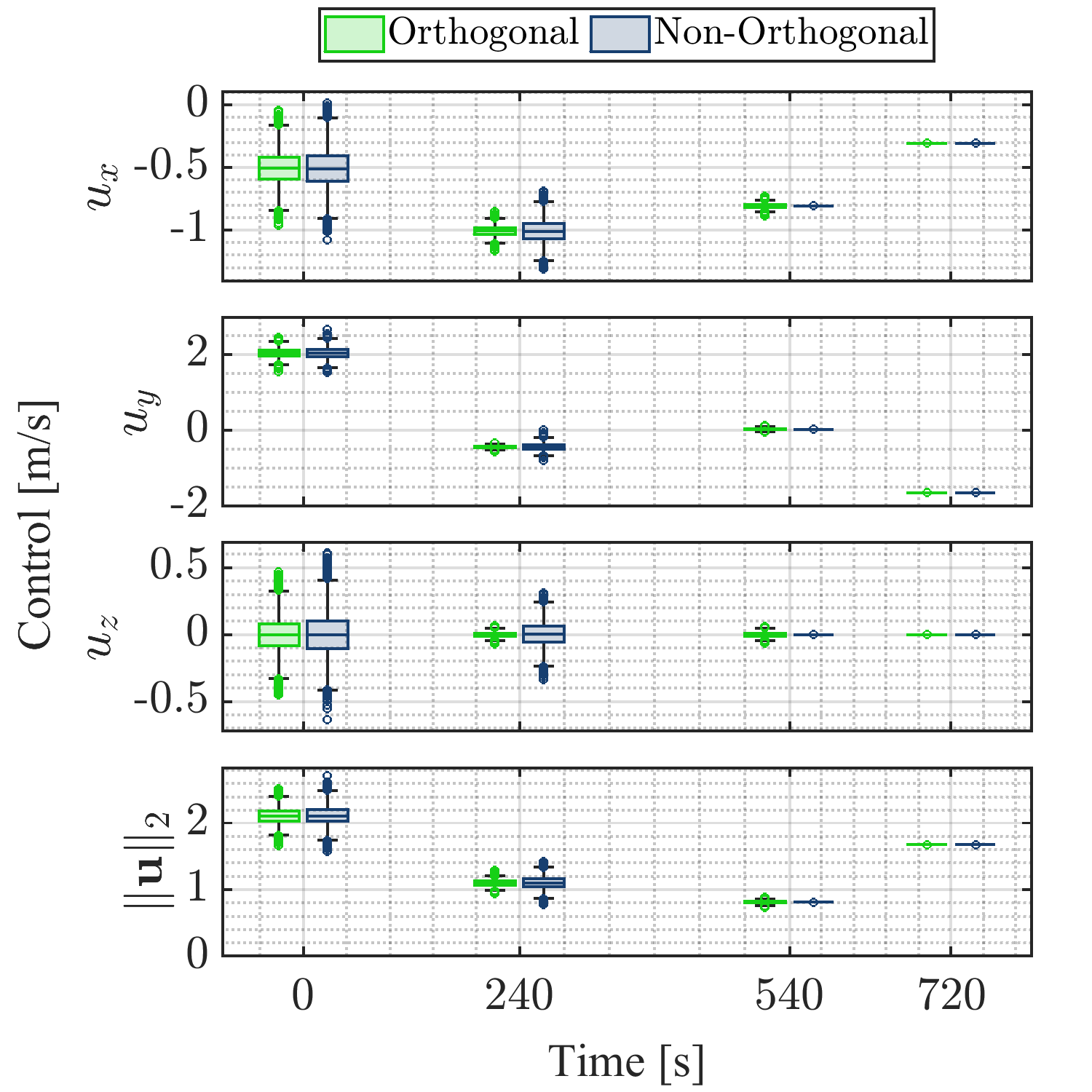}}
    \hskip 0.1truein
    \subcaptionbox
    {Histogram of Total $\Delta V$.\label{fig: control histogram}}
    {\includegraphics[width=0.48\textwidth]{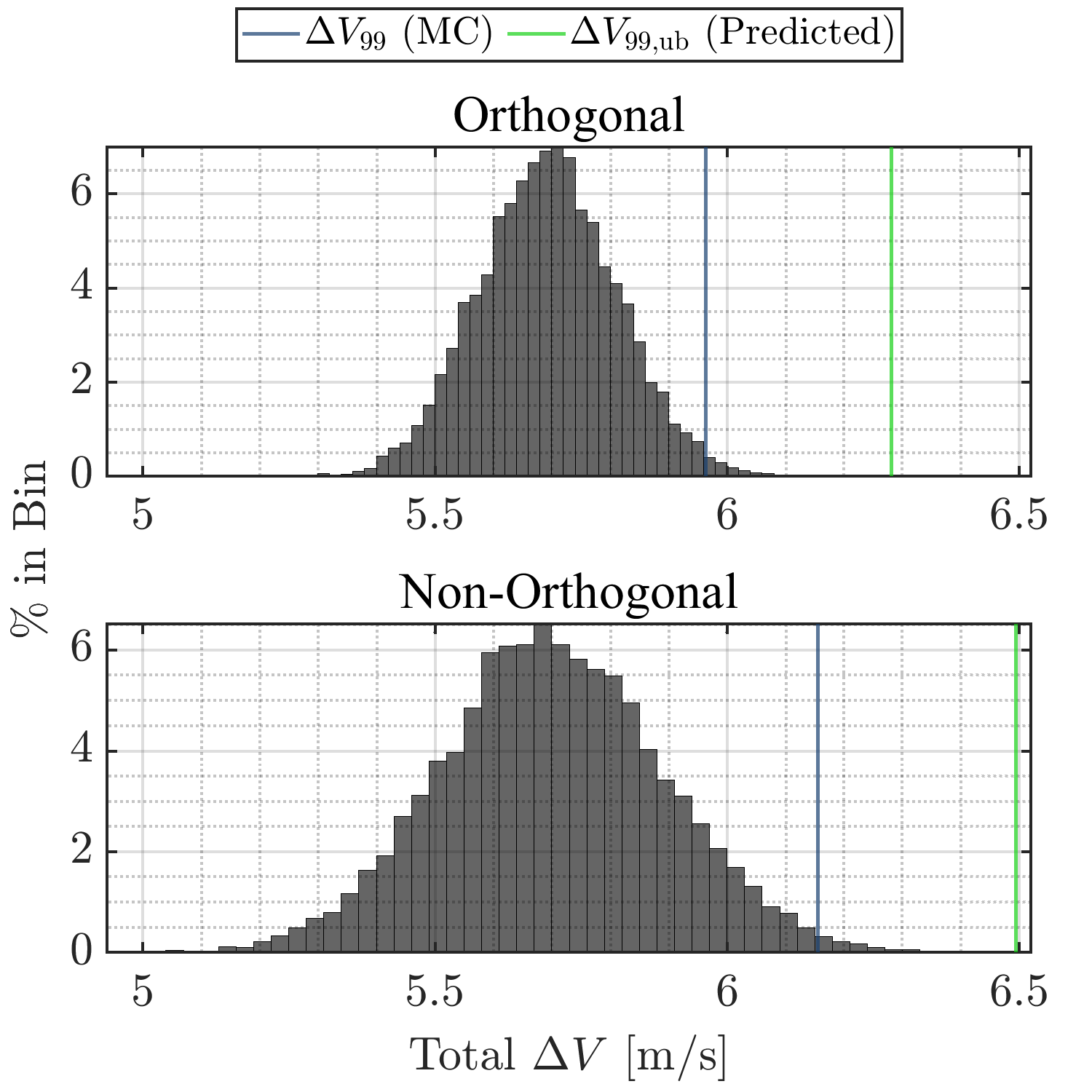}}
    \caption{Monte Carlo ($n_{\text{samples}}=10,000$): Comparison between orthogonal CS-OF and non-orthogonal CS-OF, both with measurement underweighting.}
    \label{fig: control cost}
\end{figure}

Figure~\ref{fig: control cost} shows the differences in the control statistics between the orthogonal and non-orthogonal CS-OF approaches. Due to the complete divergence of the orthogonal CS-OF without measurement underweighting case, it is omitted from this analysis. Figure~\ref{fig: control history} shows the distribution of each maneuver. The nominal maneuver is identical in both cases, but the non-orthogonal formulation appears to produce greater variance in subsequent control actions, particularly in the $x$ and $z$ components of the second maneuver. This difference likely reflects the influence of the cross-covariance terms on the control policy. Despite these differences, Figure~\ref{fig: control histogram} shows that the total $\Delta V$ costs between the two methods are largely the same. The non-orthogonal method yields a slightly higher $\Delta V_{99}$, consistent with the additional control effort required to satisfy the terminal covariance constraint, which the orthogonal method fails to meet.

\subsubsection{Filter Consistency v.s. Steering Consistency}
To further study the differences in filter consistency and steering consistency, an analysis of the underweighting factor is conducted. Both the orthogonal and non-orthogonal CS-OF problems are re-solved with underweighting factors ranging from a slightly underweighted $p = 0.95$ to an extreme $p=0.5$ case. For each underweighting case, a consistency metric is applied to the terminal true state distribution and the filtering process with full-state consideration.

\begin{figure}[htb]
	\centering\includegraphics[width=0.75\textwidth]{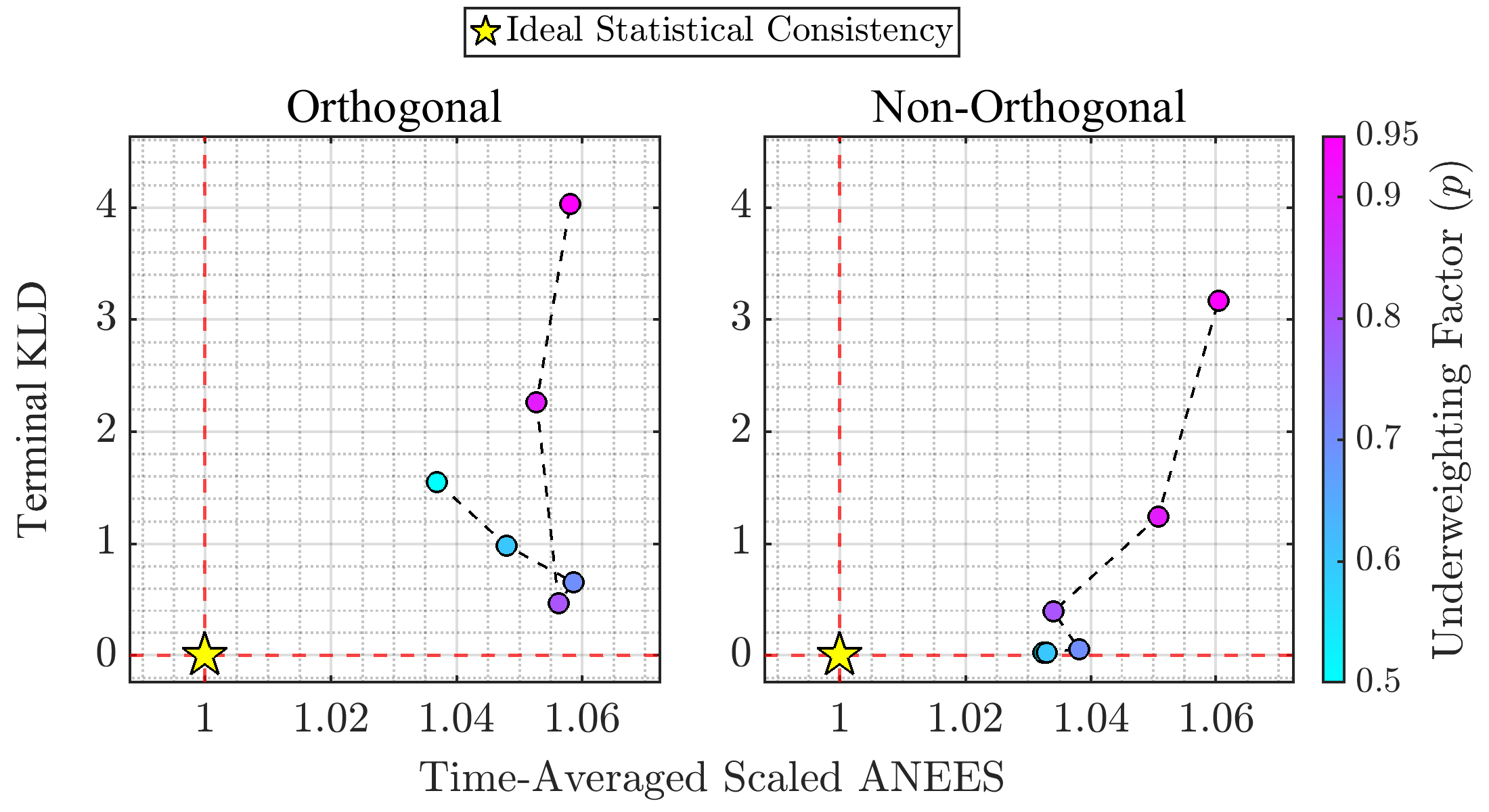}
    \caption{Monte Carlo ($n_{\text{samples}}=10,000$): Consistency metrics between orthogonal and non-orthogonal CS-OF.}
	\label{fig: consistency}
\end{figure}

The Kullback-Leibler divergence (KLD) serves as the steering consistency metric. This value is computed between the Monte Carlo and predicted distributions at terminal time, assuming Gaussian statistics. A KLD of zero represents that the two probability distributions are identical, or in other words, the actual distribution under the CS-OF policy is consistent with the one predicted by it. The Scaled Average Normalized Estimation Error Squared (ANEES) is used as the filter consistency metric \cite{Bar-Shalom-Estimation}. ANEES represents the average NEES across the Monte Carlo runs at a given time, and scaling ANEES by the state dimension $n_x$ results in a value of one for statistically consistent filters. To condense the ANEES scores into a single dimension, they are further averaged across times.\footnote{This practice of time-averaging is not recommended for rigorous filter analysis as the consistency can change through time. Nevertheless, it provides a compact summary of the overall trends in this analysis.} The relationship between the terminal KLD and time-averaged scaled ANEES scores from CS-OF solutions across different underweighting factors is shown in Figure~\ref{fig: consistency}.

From Figure~\ref{fig: consistency}, it can be seen that decreasing the underweighting factor improves the statistical consistency of the filter for both the orthogonal and non-orthogonal cases, and the improvement is about the same for each case. But the important detail is how the terminal KLD is affected by the underweighting. For orthogonal CS-OF, this decrease in underweighting initially improves the steered state's consistency, which is likely from the benefits of better filter convergence. However, there is a point where this consistency starts to depart from the ideal. Since an increase in measurement underweighting makes the filter gain increasingly suboptimal, the statistical inconsistency due to ignoring the non-orthogonal statistics starts to be more apparent in the orthogonal CS-OF solution. This is contrasted with the non-orthogonal case. Even after the initial filter improvement from underweighting is accounted for, the terminal KLD never departs from the ideal statistical consistency even with high underweighting. This behavior is due to the non-orthogonal formulation accounting for the cross-covariance terms and reducing the associated prediction errors. 

Filter consistency remains imperfect in both cases, reflecting the conservative uncertainty estimates introduced by measurement underweighting. Nevertheless, these results show that improved filter consistency does not necessarily translate into better agreement between predicted and realized state statistics, since non-orthogonality also affects covariance propagation.

\subsection{Need for Gain Correction Term}
Although the underweighted filters in both formulations converge well, discrepancies remain between the reference and realized $3\sigma$ error bounds in Figures~\ref{fig: ortho_CSOF_uw, filter} and~\ref{fig: nonortho_CSOF, filter}. These discrepancies suggest a mismatch between the reference and realized Kalman gains, motivating further analysis of gain correction. As noted earlier, gain correction is omitted in Figures~\ref{fig: ortho_CSOF} and~\ref{fig: ortho_CSOF_uw} to replicate the conventional approach, and included in Figure~\ref{fig: nonortho_CSOF} for the proposed approach. To complete the comparison, the orthogonal cases are evaluated with gain correction and the non-orthogonal case without it, all at the original underweighting factor of $p=0.7$.

\begin{figure}[htb]
	\centering\includegraphics[width=0.75\textwidth]{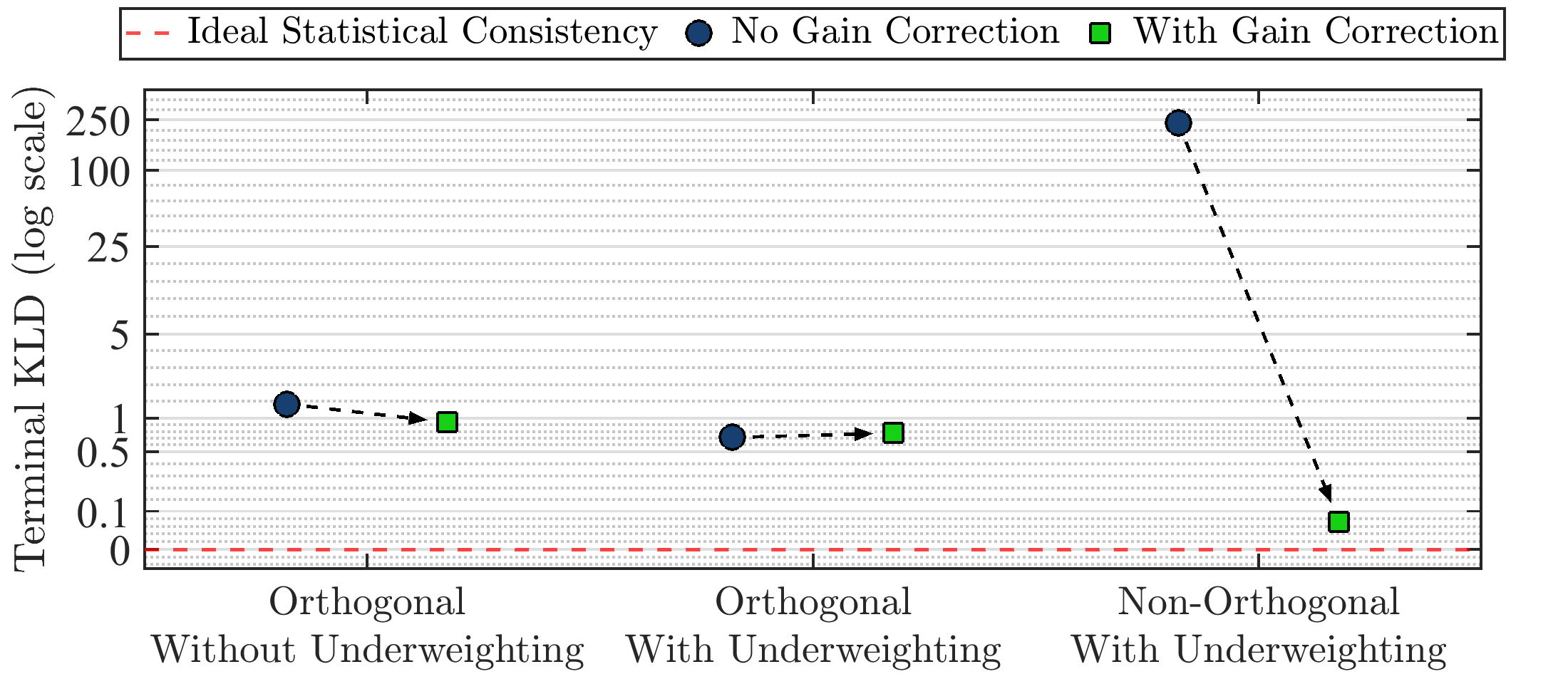}
    \caption{Monte Carlo ($n_{\text{samples}}=10,000$): Analysis of CS-OF methods with and without gain correction.}
	\label{fig: gain correction}
\end{figure}

Figure~\ref{fig: gain correction} shows the terminal KLD across all CS-OF methods, with and without the gain correction term. Firstly, it can be seen that gain correction slightly improves the statistical consistency of orthogonal CS-OF without measurement underweighting. Recalling from Figure~\ref{fig: ortho_CSOF}, many of the realized filter error is inconsistent with the reference filter error, and as such each has a different $L_k$ compared to $L_k^*$. Gain correction is meant to remedy this, hence the small improvement; but regardless, the filter still diverges which leads to a loss in controlled state accuracy. On the other hand, Figure~\ref{fig: gain correction} shows little apparent improvement from gain correction in orthogonal CS-OF with measurement underweighting. This agrees with the fact that gain correction does not alter the underlying orthogonality assumption. Together, these results suggest that gain correction alone cannot fully resolve discrepancies caused by filter divergence or an invalid orthogonality assumption.

An important result is that Figure~\ref{fig: gain correction} shows pronounced divergence in the final position distribution when gain correction is omitted from non-orthogonal CS-OF. This degree of divergence is not observed in the corresponding orthogonal CS-OF cases. In non-orthogonal CS-OF, the propagation of the cross-covariance between the state estimate and the estimation error depends on the filter gain. A mismatch between the reference and realized gains can therefore introduce errors in this term, potentially amplifying discrepancies in the predicted state covariance. Orthogonal CS-OF assumes zero cross-covariance between the estimate and the error, which may explain its lower sensitivity to filter gain mismatch in this example. This additional dependence on the filter gain may make the non-orthogonal formulation more sensitive to discrepancies between the reference and realized gains. Nevertheless, with the gain correction term, the non-orthogonal CS-OF is a magnitude better than the orthogonal cases, and almost reaches the ideal statistical consistency.  

\section{Discussion on Non-Orthogonal CS-OF}\label{sec: discussion}
Table~\ref{tab:CSOF_results} summarizes the general findings from Section~\ref{sec: numerical example}: underweighting improves filter convergence, but does not necessarily lead to accurate state predictions in the CS-OF problem. The following sections provide additional insights into the problem not covered by the numerical results.
\begin{table}[htbp]
    \fontsize{10}{10}\selectfont
    \caption{Summary of CS-OF Results from Section~\ref{sec: numerical example}}
    \label{tab:CSOF_results}
    \centering

    \newcommand{\CSOFcell}[1]{%
        \begin{tabular}[c]{@{}c@{}}
            \noalign{\vskip 3pt}
            #1 \\
            \noalign{\vskip 3pt}
        \end{tabular}%
    }

    \begin{tabular}{c | c c}
        \hline
        \CSOFcell{Measurement\\Underweighting}
            & Orthogonal & Non-Orthogonal \\
        \hline
        \CSOFcell{No Underweighting}
            & \CSOFcell{Filter diverges\\State predictions inaccurate}
            & \CSOFcell{\textemdash} \\
        \hline
        \CSOFcell{With Underweighting}
            & \CSOFcell{Filter converges\\State predictions inaccurate}
            & \CSOFcell{Statistically consistent\\filter/state predictions} \\
        \hline
    \end{tabular}
\end{table}

\subsection{Non-Orthogonal CS-OF as a Generalization of the CS-OF Problem}
The initial statistical relationships of the true and estimate are given by $\mathbf{P}_0^{-}$. As seen in the derivation, no assumption is made regarding the relationship between $\boldsymbol{x}_0$ and $\hat{\boldsymbol{x}}^-_0$ other than that they are jointly Gaussian from Assumption~\ref{assump: joint Gaussian}. This implies that the proposed non-orthogonal CS-OF framework can accommodate both orthogonal and non-orthogonal $\mathbf{P}_0^{-}$, which provides a broader class of initial conditions than conventional CS-OF formulations. The orthogonal initialization in Eq.~\eqref{eq: numerical example, initial cov} from Section~\ref{sec: numerical example} provides a common starting point for the most direct comparison with these methods; a non-orthogonal initialization would immediately fall outside orthogonal CS-OF's assumptions, even for optimal Kalman gains.

\subsection{Choice of Initial Statistical Relationship and Connection to LinCov}
The preceding section raises a conceptual question: why might the initial state estimate and its estimation error be non-orthogonal? Existing orthogonal CS-OF formulations adopt an optimal filtering framework \cite{Ridderhof-CSOF, Pilipovsky-CSOF, Naoya-CS-JGCD}, in which the state estimate and its estimation error are orthogonal. This restricts the admissible initial covariances to $P_0 \succeq \hat{P}_0$, requiring the uncertainty of the true state to be no smaller than that of the estimate. 

A related perspective arises in linear covariance analysis (LinCov), which predicts the evolution of state uncertainty under feedback control for mission analysis. Prior studies \cite{Woffinden-LinCov, Williams-LinCov, Cunningham-LinCov, Geller-LinCov} examine how measurement and estimation errors influence feedback control and the resulting state dispersion about a \emph{predetermined nominal trajectory}. When the initial estimation error is assumed to be uncorrelated with the true state, the covariance relationship becomes $\hat{P}_0 = P_0 + \widetilde{P}_0$, implying $P_0 \preceq \hat{P}_0$. Although these studies do not explicitly mention CS or CS-OF, it is clear that LinCov and covariance steering share a common focus on uncertainty propagation under control. LinCov-based optimization has been addressed using particle swarm optimization \cite{Cunningham-LinCov} and genetic algorithms \cite{Geller-LinCov}, whereas CS-OF primarily focuses on convex or sequential convex formulations. By accommodating the initial covariance relationships used in LinCov, the proposed non-orthogonal CS-OF framework offers a bridge between these approaches.

\subsection{Practical Considerations of Orthogonal and Non-Orthogonal Formulation}
Orthogonal CS-OF nevertheless retains several practical advantages. The separation of the estimation process from the optimization problem allows for a lower-dimensional optimization problem compared to non-orthogonal CS-OF. Research in orthogonal CS-OF has also developed methods to obtain Markovian feedback gains \cite{Pilipovsky-CSOF, Naoya-Markovian, Oguri-Safe-Autonomy}. Compared to the SHF policy, a Markovian policy is more computationally attractive for onboard spacecraft because it does not require onboard computers to store the past state-estimate history. A Markovian feedback synthesis method is also available for non-orthogonal CS-OF \cite{Qi-Non-Orthogonal}, although it relies on a rank-constrained sequential convex framework. Orthogonal CS-OF also draws on a substantial body of theoretical work \cite{Chen-CSOF, Liu-CS, Liu-CS-automatica}. Extending these results to non-orthogonal CS-OF is less straightforward because the cross-covariance terms introduce additional coupling between estimation and control.

Ultimately, the proposed non-orthogonal CS-OF formulation still relies on linearized models, and thus highly nonlinear cases likely exist where all CS-OF methods fail. Reducing discrepancies between the reference and realized filter gains, or reparameterizing the measurement model to reduce nonlinearity \cite{Craig-Angles}, may improve its statistical consistency and performance. The suitability of these refinements will depend on the application.

\section{Conclusion}
This paper addresses covariance steering with output-feedback (CS-OF) under nonlinear measurements. Although CS-OF is formulated for linear stochastic systems, its extension to nonlinear problems through linearization has attracted considerable interest for aerospace guidance. This paper examines two sources of statistical inconsistency in linearized CS-OF under nonlinear measurements and presents methods to address them. First, conventional CS-OF formulations do not take into account the effects of sub-optimal filter gains such as measurement underweighting. While underweighting helps avoid filter divergence under nonlinear measurements, it induces a nonzero cross-covariance between the state estimate and the estimation error, violating the orthogonality assumption underlying these formulations. This leads to the first contribution of this paper, which is a convex solution to the non-orthogonal CS-OF problem. Second, the filter gains realized during nonlinear propagation can differ from those computed within the linearized CS-OF solution. To mitigate this mismatch, the second contribution of this paper introduces a gain correction term to the optimized CS-OF policy during closed-loop operation. An angles-only spacecraft rendezvous example demonstrates that the non-orthogonal formulation and gain correction are jointly needed to obtain consistent filter statistics and predictions of state dispersion. This paper emphasizes that statistically consistent CS-OF requires thorough consideration of both the estimation and covariance steering processes, since techniques that improve the estimation may inadvertently invalidate assumptions underlying the control formulation.

\section*{Appendix}
\subsection{Sparse Structure of Block F}\label{appendix: sparse F}
\begin{lemma}\label{lemma: block F}
Given a strictly lower-triangular block matrix $\mathbf{B}$ satisfying $\mathbf{B}_{i,j} = 0$ whenever $i\leq j$ for $i=0,\ldots,N-1$ and $j=0,\ldots,N-2$. The defined decision variable $\mathbf{F} = \mathbf{K}(I-\mathbf{B}\mathbf{K})^{-1}$ follows the same sparse structure as $\mathbf{K}$:
\begin{equation}\label{eq: sparse F} 
\mathbf{K}_{k,i} =
    \begin{cases}
        \begin{bmatrix} 0_{n_u \times n_x} & K_{k,i} \end{bmatrix},
        & i \leq k \\
        0_{n_u \times 2n_x}, & i > k
    \end{cases}
\qquad\rightarrow\qquad
    \mathbf{F}_{k,i} =
    \begin{cases}
        \begin{bmatrix} 0_{n_u \times n_x} & F_{k,i} \end{bmatrix},
        & i \leq k \\
        0_{n_u \times 2n_x}, & i > k
    \end{cases}
\end{equation}
for $k=0,\ldots,N-2$ and $i=0,\ldots,N-1$.
\end{lemma}
\begin{proof}
Consider the block products
\begin{equation}
    (\mathbf{K}\mathbf{B})_{k,j}
    = \sum_{i=0}^{N-1}
      \mathbf{K}_{k,i}\mathbf{B}_{i,j}
,\qquad
(\mathbf{B}\mathbf{K})_{i,\ell}
    = \sum_{j=0}^{N-2}
      \mathbf{B}_{i,j}\mathbf{K}_{j,\ell}     
\end{equation}
For $\mathbf{K}\mathbf{B}$, a nonzero summand requires $i \leq k$ and $j<i$, or $j<k$. The zero summands are likewise when $k\leq j$, proving that $\mathbf{K}\mathbf{B}$ is strictly block lower triangular. Similarly, a nonzero summand in $\mathbf{B}\mathbf{K}$ requires $\ell \leq j<i$ with zero summands on $i \leq \ell$. Thus, $\mathbf{B}\mathbf{K}$ is also strictly block lower triangular. It follows that $I - \mathbf{K}\mathbf{B}$ and $I - \mathbf{B}\mathbf{K}$ are both invertible. Furthermore, $\mathbf{K}\mathbf{B}$ being strictly block lower triangular is nilpotent with index $N-1$. For a nilpotent matrix, the Neumann series for the matrix inverse is applied:
\begin{equation}
    \mathbf{T} = (I - \mathbf{K}\mathbf{B})^{-1} = 
    \sum_{r=0}^{N-2}(\mathbf{K}\mathbf{B})^r
\end{equation}
Thus, $\mathbf{T}$ is block lower triangular. The identity $(I-\mathbf{K}\mathbf{B})\mathbf{K} = \mathbf{K}(I-\mathbf{B}\mathbf{K})$ combined with the definition of $\mathbf{F}$ and arranging gives $\mathbf{F}= (I-\mathbf{K}\mathbf{B})^{-1}\mathbf{K} = \mathbf{T}\mathbf{K}$. Expressed as a summation and noting that $\mathbf{T}_{k,j} = 0$ for $k < j$:
\begin{equation}
\mathbf{F}_{k,i}
= \sum_{j=0}^{k}
  \mathbf{T}_{k,j}\mathbf{K}_{j,i}
\end{equation}
If $i > k$, then $i > j$ and $\mathbf{K}_{j,i} = 0$ for every term in the sum, causing the summand to be zero. If $i \leq k$, this gives
\begin{equation}
\mathbf{F}_{k,i} = 
\sum_{j=i}^{k}
   \mathbf{T}_{k,j}
   \begin{bmatrix} 0 & K_{j,i} \end{bmatrix}
   =
   \begin{bmatrix} 0 & F_{k,i} \end{bmatrix}
\end{equation}
Combining these cases yields Eq.~\eqref{eq: sparse F} and proves the claim.
\end{proof}

\subsection{Initial Augmented Covariance}\label{appendix: initial cov}
Let the augmented covariance be divided into submatrices corresponding to the covariances of the true and a priori estimate:
\begin{equation}
\mathbf{P}^-_k  =
\begin{bmatrix}
    P_k & (\Sigma^-_k)^\top \\
    \Sigma^-_k & \hat{P}^-_k
\end{bmatrix}
\end{equation}
where $\Sigma^-_k = \mathrm{Cov}(\hat{\boldsymbol{x}}^-_k,\boldsymbol{x}_k)$. At the initial time, define $\hat{P}^-_0 = \mathrm{Cov}(\hat{\boldsymbol{x}}^-_0)$ and $\widetilde{P}^-_0 = \mathrm{Cov}(\boldsymbol{x}_0 - \hat{\boldsymbol{x}}^-_0)$. If the initial estimate is unbiased and orthogonal, then $\mathrm{Cov}(\hat{\boldsymbol{x}}^-_0,\boldsymbol{x}_0 - \hat{\boldsymbol{x}}^-_0) = 0$. It follows that
\begin{equation}
    \begin{aligned}
P_0 &= \mathrm{Cov}[\hat{\boldsymbol{x}}^-_0 + (\boldsymbol{x}_0 - \hat{\boldsymbol{x}}^-_0)]=
\mathrm{Cov}(\hat{\boldsymbol{x}}^-_0) 
+
\mathrm{Cov}(\boldsymbol{x}_0 - \hat{\boldsymbol{x}}^-_0) 
+    
\mathrm{Cov}(\hat{\boldsymbol{x}}^-_0, \boldsymbol{x}_0 - \hat{\boldsymbol{x}}^-_0) 
+
\mathrm{Cov}(\hat{\boldsymbol{x}}^-_0,\boldsymbol{x}_0 - \hat{\boldsymbol{x}}^-_0)^\top = \hat{P}^-_0 + \widetilde{P}^-_0 \\
~
\Sigma^-_0 &= \mathrm{Cov}[\hat{\boldsymbol{x}}^-_0,\hat{\boldsymbol{x}}^-_0 + (\boldsymbol{x}_0 - \hat{\boldsymbol{x}}^-_0)] = 
\mathrm{Cov}(\hat{\boldsymbol{x}}^-_0,\hat{\boldsymbol{x}}^-_0)
+
\mathrm{Cov}(\hat{\boldsymbol{x}}^-_0,\boldsymbol{x}_0 - \hat{\boldsymbol{x}}^-_0) = \hat{P}^-_0
    \end{aligned}
\end{equation}
which yields the desired Eq.~\eqref{eq: numerical example, initial cov}.

\section*{Acknowledgments}
The authors would like to thank fellow lab mates Nicholas H. Craig and Naoya Kumagai for the discussions on non-orthogonality within covariance steering. This material is based on work supported by the Air Force Office of Scientific Research under award number FA9550-23-1-0512. A large language model was used solely for language editing and figure formatting. It was not used for generating content, analysis, or conclusions. All outputs were reviewed and validated by the authors.

\bibliography{references}

@article{Ridderhof-CSOF,
  author  = {Jack Ridderhof and Kazuhide Okamoto and Panagiotis Tsiotras},
  title   = {Chance Constrained Covariance Control for Linear Stochastic Systems With Output Feedback},
  journal = {IEEE Conference on Decision and Control}, 
  volume  = {},
  number = {},
  pages={1758-1763},
  year = {2020},
  doi={10.1109/CDC42340.2020.9303731}
}

@article{Qi-min-NL,
	Author = {Daniel C. Qi and Kenshiro Oguri},
	Journal = {Journal of Guidance, Control, and Dynamics},
	Pages = {3014-3024},
	Title = {Convex Bound of Nonlinear Dynamical Errors for Covariance Steering},
	Volume = {49},
	Number = {10},
	Year = {2026},
    doi = {10.2514/1.G009692}}

@article{Skaf-affine-controller,
  author  = {Skaf, Joëlle and Boyd, Stephen P.},
  title   = {Design of Affine Controllers via Convex Optimization}, 
  journal = {IEEE Transactions on Automatic Control},
  volume  = {55},
  number = {11},
  pages={2476-2487},
  year = {2010},
  doi={10.1109/TAC.2010.2046053}
}

@article{Oguri-Safe-Autonomy,
  author = {Kenshiro Oguri},
  title = {Chance-Constrained Control for Safe Spacecraft Autonomy: Convex Programming Approach},
  journal = {American Control Conference},
  volume  = {},
  number = {},
  pages={2318-2324},
  year = {2024},
  doi={10.23919/ACC60939.2024.10645008}
}

@article{Okamoto-CS,
  author={Okamoto, Kazuhide and Goldshtein, Maxim and Tsiotras, Panagiotis},
  journal={IEEE Control Systems Letters}, 
  title={Optimal Covariance Control for Stochastic Systems Under Chance Constraints}, 
  year={2018},
  volume={2},
  number={2},
  pages={266-271},
  doi={10.1109/LCSYS.2018.2826038}
}

@article{Okamoto-CS-Vehicle,
  author={Okamoto, Kazuhide and Tsiotras, Panagiotis},
  journal={IEEE Robotics and Automation Letters}, 
  title={Optimal Stochastic Vehicle Path Planning Using Covariance Steering}, 
  year={2019},
  volume={4},
  number={3},
  pages={2276-2281},
  doi={10.1109/LRA.2019.2901546}}

@article{Zanetti-Underweighting,
	Author = {Zanetti, Renato and DeMars, Kyle J. and Bishop, Robert H.},
	Journal = {Journal of Guidance, Control, and Dynamics},
	Pages = {1670-1675},
	Title = {Underweighting Nonlinear Measurements},
	Volume = {33},
	Number = {5},
	Year = {2010},
    doi = {10.2514/1.50596}}

@article{Woffinden-LinCov,
  author  = {David C. Woffinden},
  title   = {Linear Covariance Techniques To Analyze A Multi-vehicle, Multi-GN\&C System With Applications To Rendezvous In A Near Rectilinear Halo Orbit},
  journal = {AAS Guidance, Navigation and Control Conference},
  Pages = {1-24},
  volume  = {},
  number = {},
  year = {2023},
  doi = {}
}

@article{Cunningham-LinCov,
  author  = {David A. Cunningham and Ryan P. Russell and David C. Woffinden},
  title   = {Robust Trajectory Optimization for NRHO Rendezvous Using SPICE Kernel Relative Motion},
  journal = {AAS Guidance, Navigation and Control Conference},
  Pages = {1-17},
  volume  = {},
  number = {},
  year = {2025},
  doi = {}
}

@article{Williams-LinCov,
  author  = {James Williams and William E. Brandenburg and David Woffinden and Zachary R. Putnam},
  title   = {Validation of Linear Covariance Techniques for Mars Entry, Descent, and Landing Guidance and Navigation Performance Analysis},
  journal = {AIAA SciTech},
  Pages = {},
  volume  = {},
  number = {},
  year = {2021},
  doi = {10.2514/6.2022-0745}
}

@book{Maybeck-StochasticTextbook,
  author    = {Peter S. Maybeck},
  title     = {Stochastic Models, Estimation, and Control},
  publisher = {Academic Press},
  address   = {New York},
  year      = {1979}
}

@article{Naoya-CS-JGCD,
  author  = {Naoya Kumagai and Kenshiro Oguri},
  title   = {Robust Cislunar Low-Thrust Trajectory Optimization under Uncertainties via Sequential Covariance Steering},
  journal = {Journal of Guidance, Control, and Dynamics},
  pages = {2725-2743},
  volume  = {48},
  number = {12},
  year = {2025},
  doi = {10.2514/1.G009092}
}

@article{Pilipovsky-CSOF,
  author  = {Pilipovsky, Joshua and Tsiotras, Panagiotis},
  title   = {Computationally Efficient Chance Constrained Covariance Control with Output Feedback}, 
  journal = {IEEE Conference on Decision and Control}, 
  volume  = {},
  number = {},
  pages= {677-682},
  year = {2024},
  doi={10.1109/CDC56724.2024.10885876}
}

@article{Hotz-CS,
  author  = {Anthony Hotz and Robert E. Skelton},
  title   = {Covariance control theory},
  journal = {International Journal of Control},
  volume  = {46},
  number = {1},
  pages= {13--32},
  year = {1987},
  doi={10.1080/00207178708933880}
}

@article{Chen-CS,
  author  = {Yongxin Chen and Tryphon T. Georgiou and Michele Pavon},
  title   = {Optimal Steering of a Linear Stochastic System to a Final Probability Distribution, Part 1},
  journal = {IEEE Transactions on Automatic Control}, 
  volume  = {61},
  number = {5},
  pages={1158-1169},
  year = {2016},
  doi={10.1109/TAC.2015.2457784}
}

@article{Ridderhof-SCP,
	Author = {Jack Ridderhof and Joshua Pilipovsky and Panagiotis Tsiotras},
    Title = {Chance-Constrained Covariance Control for Low-thrust Minimum-Fuel Trajectory Optimization},
	Journal = {AIAA/AAS Astrodynamics Specialists Conference},
    volume  = {},
    number = {},
	Pages = {1--20},
	Year = {2021},
    doi = {}}

@article{Jerry-CovSteering,
	Author = {Jerry Varghese and Kenshiro Oguri},
	Journal = {Journal of Guidance, Control, and Dynamics},
    Pages = {59-77},
	Title = {Trajectory Optimization Under Uncertainty with Nonlinear Programming and Forward–Backward Shooting},
	Volume = {49},
	Number = {1},
	Year = {2026},
    doi = {10.2514/1.G009259}}

@article{Benedikter-CovSteering,
    Author = {Boris Benedikter and Alessandro Zavoli and Zhenbo Wang and Simone Pizzurro and Enrico Cavallini},
    Journal = {Journal of Guidance, Control, and Dynamics},
    Pages = {2061--2075},
    Title = {Convex Approach to Covariance Control with Application to Stochastic Low-Thrust Trajectory Optimization},
	Volume = {45},
	Number = {11},
	Year = {2022},
    doi = {10.2514/1.G006806}}

@article{Fife-CS,
  author  = {William N Fife and Kyle J. DeMars},
  journal = {The Journal of the Astronautical Sciences},
  title  = {A Minimum Initial Information Approach for Nominal Guidance via Convex Optimization},
  Pages = {1-27},
  volume  = {72},
  number = {18},
  year = {2025},
  doi = {10.1007/s40295-025-00494-7}
}

@article{Alkhenani-CS,
	Author = {Abdullah K. Alkhenani and Kenshiro Oguri},
	Journal = {Journal of Guidance, Control, and Dynamics},
    Pages = {1-18},
	Title = {Chance-Constrained Autonomous Spacecraft Formation Reconfiguration and Rendezvous with Passive Abort Safety},
	Volume = {},
	Number = {},
	Year = {2026},
    doi = {10.2514/1.G010121}}

@article{Xiao-CS,
	Author = {Xiao, Yizheng and Gong, Youmin and Mei, Jie and Ma, Guangfu and Wu, Weiren},
	Journal = {Journal of Guidance, Control, and Dynamics},
    Pages = {1-18},
	Title = {Covariance Steering for Powered-Descent Guidance with Embedded Navigation and Tracking},
	Volume = {},
	Number = {},
	Year = {2026},
    doi = {10.2514/1.G010096}}

@article{Qi-CUT-Steering,
	Author = {Daniel C. Qi and Kenshiro Oguri and Puneet Singla and Maruthi R. Akella},
    title  = {Non-Gaussian Distribution Steering in Nonlinear Dynamics with Conjugate Unscented Transformation},
	Journal = {Journal of Guidance, Control, and Dynamics},
    Pages = {},
	Volume = {},
	Number = {},
	Year = {in press},
    doi = {}}

@article{Burnett-NG-Control,
	Author = {Burnett, Ethan and Michelotti, Niccolò and Topputo, Francesco},
    title  = {Confidence Boundary Methods for Non-Gaussian Mission Analysis and Stochastic Guidance},
	Journal = {ESA GNC \& ICATT Conference},
    Pages = {},
	Volume = {},
	Number = {},
	Year = {2026},
    doi = {}}

@article{Shuttle-Underweight,
	Author = {Kriegsman, Bernard A. and Tao, Yee-Chee},
    title  = {Shuttle Navigation System for Entry and Landing Mission Phases},
	Journal = {Journal of Spacecraft and Rockets},
    pages = {213-219},
    volume = {12},
    number = {4},
	Year = {1975},
    doi = {10.2514/3.56966}}

@article{Park-Underweight,
	Author = {Park, Young Bum and Park, Chan Gook},
    title  = {Mitigation of Vision Measurement Nonlinearity Effect on Lunar Descent Navigation Using Underweighting},
	Journal = {Journal of Guidance, Control, and Dynamics},
    pages = {2370-2377},
    volume = {40},
    number = {9},
	Year = {2017},
    doi = {10.2514/1.G002729}}

@article{Liao-Underweight,
	Author = {Liao, Aaron R. and Oguri, Kenshiro},
    title  = {Autonomous Shape-Free Small-Body Navigation Using Relative Angles-Only Measurements},
	Journal = {Journal of Spacecraft and Rockets},
    pages = {677-685},
    volume = {63},
    number = {2},
	Year = {2026},
    doi = {10.2514/1.A36415}}

@article{Naoya-Sqrt-CS,
  author={Kumagai, Naoya and Oguri, Kenshiro},
  journal={IEEE Transactions on Automatic Control}, 
  title={Square Root-Factorized Covariance Steering}, 
  year={2026},
  volume={},
  number={},
  pages={1-8},
  doi={10.1109/TAC.2026.3738329}}

@article{Chen-CSOF,
  author={Chen, Yongxin and Georgiou, Tryphon and Pavon, Michele},
  journal={IEEE Conference on Decision and Control}, 
  title={Steering state statistics with output feedback}, 
  year={2015},
  volume={},
  number={},
  pages={6502-6507},
  doi={10.1109/CDC.2015.7403244}}

@article{Rapakoulias-CS-CC,
  author={Rapakoulias, George and Tsiotras, Panagiotis},
  journal={IEEE Conference on Decision and Control}, 
  title={Discrete-Time Optimal Covariance Steering via Semidefinite Programming}, 
  year={2023},
  volume={},
  number={},
  pages={1802-1807},
  doi={10.1109/CDC49753.2023.10384118}}

@article{Geller-LinCov,
  author  = {David Geller and Simon Shuster and David Woffinden and Stefan Bieniawski},
  title   = {Robust Cislunar Trajectory Optimization Via Midcourse Correction and Optical Navigation Scheduling},
  journal = {AAS Guidance, Navigation and Control Conference},
  Pages = {1-13},
  volume  = {},
  number = {},
  year = {2022},
  doi = {}
}

@article{Qi-Non-Orthogonal,
  author={Daniel C Qi and Kenshiro Oguri},
  journal={arXiv preprint}, 
  title={Generalizing Output-Feedback Covariance Steering to Incorporate Non-Orthogonal Estimation Errors}, 
  year={2026},
  volume={},
  number={},
  pages={1-7},
  doi={10.48550/arXiv.2603.29753}}

@article{Naoya-Markovian,
  author={Naoya Kumagai, Kenshiro Oguri},
  journal={arXiv preprint}, 
  title={On the Optimality of Markovian Policies for Chance-Constrained Covariance Steering}, 
  year={2026},
  volume={},
  number={},
  pages={1-12},
  doi={10.48550/arXiv.2608.22589}}

@article{Liu-CS,
  author={Liu, Fengjiao and Rapakoulias, George and Tsiotras, Panagiotis},
  journal={IEEE Transactions on Automatic Control}, 
  title={Optimal Covariance Steering for Discrete-Time Linear Stochastic Systems}, 
  year={2025},
  volume={70},
  number={4},
  pages={2289-2304},
  doi={10.1109/TAC.2024.3472788}}

@article{Liu-CS-automatica,
author = {Fengjiao Liu and Panagiotis Tsiotras},
journal = {Automatica},
title = {Reachability and controllability analysis of the state covariance for linear stochastic systems},
volume = {189},
number = {},
pages = {113036},
year = {2026},
doi = {10.1016/j.automatica.2026.113036},
}

@article{Craig-Angles,
author = {Craig, Nicholas H. and Oguri, Kenshiro},
title = {Analytical Reformulation of Angles-Only Measurements with Application to Cislunar Object Tracking},
journal = {Journal of Spacecraft and Rockets},
volume = {63},
number = {3},
pages = {916-932},
year = {2026},
doi = {10.2514/1.A36471}
}

@book{Stat-OD,
  author    = {Byron D. Tapley and Bob E. Schutz and George H. Born},
  title     = {Statistical Orbit Determination},
  publisher = {Academic Press},
  address   = {Boston},
  year      = {2004}
}

@article{Kruger-angles,
author = {Justin Kruger and Simone D’Amico},
title = {Autonomous angles-only multitarget tracking for spacecraft swarms},
journal = {Acta Astronautica},
volume = {189},
number = {},
pages = {514-529},
year = {2021},
doi = {10.1016/j.actaastro.2021.08.049},
}

@book{Bar-Shalom-Estimation,
  author    = {Yaakov Bar-Shalom and X.-Rong Li and Thiagalingam Kirubarajan},
  title     = {Estimation with Applications to Navigation and Tracking},
  publisher = {Wiley},
  address   = {New York},
  year      = {2001 }
}

@article{Luca-CS,
author = {Luca Frassinella and Kenshiro Oguri and Francesco Topputo},
title = {Safe Autonomous Guidance Under Uncertainty via Sequential Covariance Steering for Small-Body Proximity Operations},
journal = {International Symposium on Space Flight Dynamics},
volume = {},
number = {},
pages = {1-10},
year = {2026},
doi = {},
}

@article{Gabriele-CS,
author = {Gabriele Nuccio and Kenshiro Oguri and Francesco Topputo},
title = {Statistical Trajectory Optimization Under Uncertainty for Icy Moon Exploration with Probability of Impact Constraints},
journal = {International Symposium on Space Flight Dynamics},
volume = {},
number = {},
pages = {1-11},
year = {2026},
doi = {},
}

@article{Qi-Risk-Measure,
author = {Daniel Qi and Kenshiro Oguri},
title = {Safe Maneuver Planning With Cubature-based Risk Measures of Non-Gaussian Distributions in Nonlinear Systems},
journal = {AAS Guidance, Navigation and Control Conference},
volume = {},
number = {},
pages = {1-22},
year = {2026},
doi = {},
}

\end{document}